\documentclass[11pt]{article}
\usepackage{mathrsfs}
\usepackage{amssymb}
\usepackage{amsmath}
\usepackage[all]{xy}
\def\Im{\mathop{\rm Im}\nolimits}
\def\Ker{\mathop{\rm Ker}\nolimits}
\def\Coker{\mathop{\rm Coker}\nolimits}

\def\mod{\mathop{\rm mod}\nolimits}
\def\Mod{\mathop{\rm Mod}\nolimits}

\def\sup{\mathop{\rm sup}\nolimits}

\def\add{\mathop{\rm add}\nolimits}
\def\Add{\mathop{\rm Add}\nolimits}

\def\dim{\mathop{\rm dim}\nolimits}
\def\codim{\mathop{\rm codim}\nolimits}
\def\Hom{\mathop{\rm Hom}\nolimits}
\def\Ext{\mathop{\rm Ext}\nolimits}
\def\sup{\mathop{\rm sup}\nolimits}
\def\lim{\mathop{\underrightarrow{\rm lim}}\nolimits}

\def\res{\mathop{\rm res}\nolimits}
\def\cores{\mathop{\rm cores}\nolimits}
\def\gen{\mathop{\rm gen}\nolimits}
\def\cogen{\mathop{\rm cogen}\nolimits}

\title{\Large \bf Proper Resolutions and Gorenstein Categories
}
\author{Zhaoyong Huang\thanks{E-mail address: huangzy@nju.edu.cn}\\
{\it \footnotesize Department of Mathematics, Nanjing University,
Nanjing 210093, Jiangsu Province, P. R. China}}
\date{ }
\begin{document}
\baselineskip=18pt \maketitle

\begin{abstract}
Let $\mathscr{A}$ be an abelian category and $\mathscr{C}$ an
additive full subcategory of $\mathscr{A}$. We provide a method to
construct a proper $\mathscr{C}$-resolution (resp. coproper
$\mathscr{C}$-coresolution) of one term in a short exact sequence in
$\mathscr{A}$ from that of the other two terms. By using these
constructions, we answer affirmatively an open question on the
stability of the Gorenstein category $\mathcal{G}(\mathscr{C})$
posed by Sather-Wagstaff, Sharif and White; and also prove that
$\mathcal{G}(\mathscr{C})$ is closed under direct summands. In
addition, we obtain some criteria for computing the
$\mathscr{C}$-dimension and the $\mathcal{G}(\mathscr{C)}$-dimension
of an object in $\mathscr{A}$.
\end{abstract}

\vspace{0.5cm}

\centerline{\large \bf 1. Introduction}

\vspace{0.2cm}

Auslander and Bridger generalized in [AB] finitely generated
projective modules to finitely generated modules of G-dimension zero
over (commutative) Noetherian rings. Furthermore, Enochs and Jenda introduced in
[EJ1] Gorenstein projective modules for arbitrary modules over a
general ring, which is a generalization of finitely generated modules of G-dimension
zero. Also in [EJ1], Gorenstein injective modules were introduced as
the dual of Gorenstein projective modules. It is well known that
the class of modules of G-dimension zero and that of Gorenstein projective modules
coincide for finitely generated modules
over a left and right Noetherian ring, and that Gorenstein projective
modules and Gorenstein injective modules share many nice properties
of projective modules and injective modules, respectively (cf. [EJ1,
EJ2, H]). The homological properties of Gorenstein projective
and injective modules and some related generalized versions have been studied by
many authors, see [AB, AM, CFH, CFrH, CI, EJ1, EJ2, EJL, GD, H, HJ, HW,
LHX, SSW1, SSW2, TW, W], and the literatures listed in them.

Let $\mathscr{A}$ be an abelian category and $\mathscr{C}$ an
additive full subcategory of $\mathscr{A}$. Sather-Wagstaff, Sharif
and White introduce the Gorenstein category
$\mathcal{G}(\mathscr{C})$, which unifies the following notions:
modules of Gorenstein dimension zero ([AB]), Gorenstein projective
modules, Gorenstein injective modules ([EJ1]), $V$-Gorenstein projective
modules, $V$-Gorenstein injective modules ([EJL]), and so on. The {\it
Gorenstein subcategory} $\mathcal{G}(\mathscr{C})$ of $\mathscr{A}$
is defined as $\mathcal{G}(\mathscr{C})=\{M$ is an object in
$\mathscr{A}\mid$ there exists an exact sequence $\cdots \to C_1 \to
C_0 \to C^0 \to C^1 \to \cdots$ in $\mathscr{C}$, which is both
$\Hom_{\mathscr{A}}(\mathscr{C},-)$-exact and
$\Hom_{\mathscr{A}}(-,\mathscr{C})$-exact, such that $M\cong
\Im(C_0\to C^0)\}$. Set $\mathcal{G}^0(\mathscr{C})=\mathscr{C}$,
$\mathcal{G}^1(\mathscr{C})=\mathcal{G}(\mathscr{C})$, and
inductively set
$\mathcal{G}^{n+1}(\mathscr{C})=\mathcal{G}(\mathcal{G}^n(\mathscr{C}))$
for any $n\geq 1$. They proved that when $\mathscr{C}$ is
self-orthogonal, $\mathcal{G}(\mathscr{C})$ possesses many nice
properties. For example, in this case, $\mathcal{G}(\mathscr{C})$ is
closed under extensions and $\mathscr{C}$ is a projective generator
and an injective cogenerator for $\mathcal{G}(\mathscr{C})$, which
induce that $\mathcal{G}^n(\mathscr{C})=\mathcal{G}(\mathscr{C})$
for any $n\geq 1$. Also in this case, they proved that
$\mathcal{G}(\mathscr{C})$ is closed under direct summands. In
particular, they posed the following open question:

\vspace{0.2cm}

{\bf Question 1.1.} (see [SSW1, Question 5.8]) Must there be an
equality $\mathcal{G}^n(\mathscr{C})=\mathcal{G}(\mathscr{C})$ for
any $n\geq 1$?

\vspace{0.2cm}

In this paper, we will prove that the answer to this question is
affirmative. The used methods are to construct a certain proper
resolution (resp. coproper coresolution) of one term in a short
exact sequence from that of the other two terms. This paper is
organized as follows.

In Section 2, we give some terminology and some preliminary results.

In Section 3, we provide a method to construct a proper resolution
(resp. coproper coresolution) of the first or last term in a short
exact sequence from that of the other two terms. We will prove the
following two theorems and their dual results.

\vspace{0.2cm}

{\bf Theorem 1.2.} {\it Let $\mathscr{C}$ be a full subcategory of
an abelian category $\mathscr{A}$ and $0\to X \to X^0 \to X^1 \to 0$
an exact sequence in $\mathscr{A}$. Let
$$\cdots \to C_i^j \to \cdots \to C_1^j \to C_0^j \to X^j \to 0 \eqno{(1.1(j))}$$
be a proper $\mathscr{C}$-resolution of $X^j$ for $j=0,1$. Then

(1) We get the following exact sequences:
$$\cdots \to C_{i+1}^1\bigoplus C_i^0 \to \cdots
\to C_2^1\bigoplus C_1^0 \to C \to X \to 0 \eqno{(1.2)}$$ and
$$0\to C \to C_1^1\bigoplus C_0^0 \to C_0^1 \to 0$$

(2) If the exact sequence $(1.1(j))$ is
$\Hom_{\mathscr{A}}(-,\mathscr{C})$-exact for $j=0,1$, then so is
the exact sequence (1.2).

(3) If $\mathscr{C}$ is closed under finite direct sums and kernels
of epimorphisms, then the exact sequence (1.2) is a proper
$\mathscr{C}$-resolution of $X$.}

\vspace{0.2cm}

{\bf Theorem 1.3.} {\it Let $\mathscr{C}$ be a full subcategory of
an abelian category $\mathscr{A}$ and
$$0\to X_1 \to X_0 \to X \to 0 \eqno{(1.3)}$$ an exact sequence
in $\mathscr{A}$. Let
$$C_j^n \to \cdots \to C_j^1 \to C_j^0 \to X_j \to 0 \eqno{(1.4(j))}$$
be a proper $\mathscr{C}$-resolution of $X_j$ for $j=0,1$. Then

(1) We get the following exact sequence:
$$C_0^n\bigoplus C_1^{n-1} \to \cdots
\to C_0^2\bigoplus C_1^1 \to C_0^1\bigoplus C_1^0 \to C_0^0 \to X
\to 0 \eqno{(1.5)}$$

(2) If all the exact sequences (1.3) and $(1.4(j))$ are
$\Hom_{\mathscr{A}}(-,\mathscr{C})$-exact for $j=1,2$, then so is
the exact sequence (1.5).

(3) If the exact sequence (1.3) is
$\Hom_{\mathscr{A}}(\mathscr{C},-)$-exact, then the exact sequence
(1.5) is a proper $\mathscr{C}$-resolution of $X$.}

\vspace{0.2cm}

As applications of these constructions, in Section 4 we prove the following
result, in which the first assertion answers Question 1.1
affirmatively.

\vspace{0.2cm}

{\bf Theorem 1.4.} {\it Let $\mathscr{A}$ be an abelian category and
$\mathscr{C}$ an additive full subcategory of $\mathscr{A}$. Then we
have

(1) $\mathcal{G}^n(\mathscr{C})=\mathcal{G}(\mathscr{C})$ for any
$n\geq 1$.

(2) $\mathcal{G}(\mathscr{C})$ is closed under direct summands.}

\vspace{0.2cm}

Let $\mathscr{A}$ be an abelian category and $\mathscr{C}$ a full
subcategory of $\mathscr{A}$. For a positive integer $n$, an object
$A$ in $\mathscr{A}$ is called an {\it $n$-$\mathscr{C}$-syzygy
object} (of an object $M$) if there exists an exact sequence $0\to A
\to C_{n-1} \to \cdots \to C_1 \to C_0\to M \to 0$ in $\mathscr{A}$
with all $C_i$ objects in $\mathscr{C}$. For an object $M$ in
$\mathscr{A}$, $\mathscr{C}$-$\dim M$ is defined as inf$\{n\geq
0\mid$ there exists an exact sequence $0 \to C_{n} \to \cdots \to
C_{1} \to C_{0} \to M \to 0$ in $\mathscr{A}$ with all $C_i$ objects
in $\mathscr{C}\}$. We set $\mathscr{C}$-$\dim M$ infinity if no
such integer exists.

In Section 5, assume that $\mathscr{C}$ is closed under extensions.
We first prove that for a positive integer $n$, an object $A$ in
$\mathscr{A}$ is $n$-$\mathscr{C}$-syzygy if and only if it is
$n$-$\cogen\mathscr{C}$-syzygy, where $\cogen\mathscr{C}$ is a
cogenerator for $\mathscr{C}$. Next we prove that if $\mathscr{X}$
is a generator-cogenerator for $\mathscr{C}$, then for any object
$M$ in $\mathscr{A}$ and $n\geq 0$, $\mathscr{C}$-$\dim M \leq n$ if
and only if for every non-negative integer $t$ such that $0 \leq t
\leq n$, there exists an exact sequence: $0\to X_{n}\to \cdots \to
X_1\to X_0\to M\to 0$ in $\mathscr{A}$ such that $X_t$ is an object
in $\mathscr{C}$ and all $X_i$ are objects in $\mathscr{X}$ for $i
\neq t$. As a consequence, when $\mathscr{C}$ is self-orthogonal, we
obtain some criteria for computing $\mathcal{G}(\mathscr{C)}$-$\dim
M$ if it is finite.

\vspace{0.5cm}

\centerline{\large \bf 2. Preliminaries}

\vspace{0.2cm}

Throughout this paper, $\mathscr{A}$ is an abelian category, all
subcategories are full subcategories of $\mathscr{A}$ closed under
isomorphisms. We fix a subcategory $\mathscr{C}$ of $\mathscr{A}$.

In this section, we give some terminology and some preliminary
results.

\vspace{0.2cm}

{\bf Definition 2.1.} ([E]) Let $\mathscr{C}\subseteq\mathscr{D}$ be
subcategories of $\mathscr{A}$. The morphism $f: C\to D$ in
$\mathscr{A}$ with $C$ an object in $\mathscr{C}$ and $D$ an object
in $\mathscr{D}$ is called a {\it $\mathscr{C}$-precover} of $D$
if for any morphism $g: C^{'} \to D$ in $\mathscr{A}$ with $C^{'}$
an object in $\mathscr{C}$, there exists a morphism $h: C^{'}\to C$
such that the following diagram commutes:
$$\xymatrix{ & C^{'} \ar[d]^{g} \ar@{-->}[ld]_{h}\\
C \ar[r]^{f} & D}$$ The morphism $f: C\to D$ is called {\it
right minimal} if a morphism $h: C\to C$ is an automorphism whenever
$f=fh$. A $\mathscr{C}$-precover $f: C\to D$ is called a {\it
$\mathscr{C}$-cover} if $f$ is right minimal. Dually, the notions of
a {\it $\mathscr{C}$-preenvelope}, a {\it left minimal morphism} and
{\it a $\mathscr{C}$-envelope} are defined. Following Auslander and
Reiten's terminology in [AR1], a $\mathscr{C}$-(pre)cover and a
$\mathscr{C}$-(pre)envelope are called a {\it (minimal) right
$\mathscr{C}$-approximation} and a {\it (minimal) left
$\mathscr{C}$-approximation}, respectively.

\vspace{0.2cm}

Recall that an exact sequence in
$\mathscr{A}$ is called {\it $\Hom_{\mathscr{A}}(\mathscr{C},
-)$-exact} if it remains still exact after
applying the functor $\Hom_{\mathscr{A}}(\mathscr{C},-)$. Let $M$ be an object in $\mathscr{A}$.
An exact sequence:
$$\cdots \buildrel {f_{i+1}} \over \longrightarrow C_i \buildrel {f_{i}} \over \longrightarrow
\cdots \buildrel {f_{2}} \over \longrightarrow C_1 \buildrel {f_{1}} \over \longrightarrow C_0
\buildrel {f_{0}} \over \longrightarrow M \to 0$$
in $\mathscr{A}$ with all $C_i$ objects in $\mathscr{C}$ is called a
{\it $\mathscr{C}$-resolution} of $M$. Recall from [AM] that the above exact sequence is called
a {\it proper $\mathscr{C}$-resolution} of $M$ if it is a $\mathscr{C}$-resolution
of $M$ and is $\Hom_{\mathscr{A}}(\mathscr{C}, -)$-exact, that is, each $f_i$
is an epic $\mathscr{C}$-precover of $\Im f_i$; and it is called a {\it minimal proper
$\mathscr{C}$-resolution} of $M$ if each $f_i$
is an epic $\mathscr{C}$-cover of $\Im f_i$. Dually,
the notions of a {\it $\Hom_{\mathscr{A}}(-, \mathscr{C})$-exact
exact sequence}, a {\it $\mathscr{C}$-coresolution} and a {\it
(minimal) coproper $\mathscr{C}$-coresolution} of $M$ are defined.

We now introduce the notion of strongly (co)proper (co)resolutions as follows.

\vspace{0.2cm}

{\bf Definition 2.2.} Let $M$ be an object in $\mathscr{A}$.

(1) A sequence:
$$\cdots \to X_i \to \cdots \to X_1 \to X_0\to M\to 0$$ in $\mathscr{A}$ is called
{\it strongly $\Hom_{\mathscr{A}}(\mathscr{C},-)$-exact exact} if it
is exact and $\Ext_{\mathscr{A}}^1(\mathscr{C},K_i)=0$ for any
$i\geq 1$, where $K_i=\Im (X_i \to X_{i-1})$. Dually, the notion of
a {\it strongly $\Hom_{\mathscr{A}}(-,\mathscr{C})$-exact exact
sequence} is defined.

(2) An exact sequence:
$$\cdots \to C_i \to \cdots \to C_1 \to C_0 \to M \to 0$$ in $\mathscr{A}$
is called a {\it strongly proper $\mathscr{C}$-resolution} of $M$ if
it is a $\mathscr{C}$-resolution of $M$ and is strongly
$\Hom_{\mathscr{A}}(\mathscr{C},-)$-exact. Dually, the notion of a
{\it strongly coproper $\mathscr{C}$-coresolution} of $M$ is
defined.

\vspace{0.2cm}

{\it Remark 2.3.} (1) It is easy to see that a strongly (co)proper
$\mathscr{C}$-(co)resolution is a (co)proper
$\mathscr{C}$-(co)resolution. But the converse does not hold true in
general. For example, let $\mathscr{C}$ be a full subcategory of
$\mathscr{A}$ closed under finite direct sums such that there exists
an object $M$ in $\mathscr{C}$ with $\Ext_{\mathscr{A}}^1(M, M)\neq
0$. Then the exact sequence:
$$0\to M  \buildrel {\binom {1_M} 0} \over \longrightarrow
M\bigoplus M \buildrel {(0, 1_M)} \over \longrightarrow M \to 0$$ is
both a proper $\mathscr{C}$-resolution and a coproper
$\mathscr{C}$-coresolution of $M$, but it neither a strongly proper
$\mathscr{C}$-resolution nor a strongly coproper
$\mathscr{C}$-coresolution of $M$.

(2) For a ring $R$, we use $\Mod R$ to denote the category of left
$R$-modules. We have that any projective resolution (resp. injective
coresolution) of a left $R$-module $M$ is just a strongly proper
$\mathcal{P}(\Mod R)$-resolution (resp. strongly coproper
$\mathcal{I}(\Mod R)$-coresolution) of $M$, where
$\mathcal{P}(\Mod R)$ (resp. $\mathcal{I}(\Mod R)$) is the
subcategory of $\Mod R$ consisting of projective (resp. injective)
modules.

(3) By the Wakamatsu's lemma (see [X, Lemma 2.1.1]), if
$\mathscr{C}$ is closed under extensions, then a minimal
(co)proper $\mathscr{C}$-resolution of an object $M$ is a strongly
(co)proper $\mathscr{C}$-resolution of $M$.

\vspace{0.2cm}

The following two observations are useful in next section.

\vspace{0.2cm}

{\bf Lemma 2.4.} {\it Let
$$\xymatrix{M \ar[d]_{g_1}\ar[r]^{f_1} & N\ar[d]_{g}\\
X \ar[r]^f & Y}$$ be a commutative diagram in $\mathscr{A}$ and $C$
an object in $\mathscr{A}$.

(1) If this diagram is a pull-back diagram of $f$ and $g$ and
$\Hom_{\mathscr{A}}(C, g)$ is epic, then $\Hom_{\mathscr{A}}(C,
g_1)$ is also epic.

(2) If this diagram is a push-out diagram of $f_1$ and $g_1$ and
$\Hom_{\mathscr{A}}(g_1, C)$ is epic, then $\Hom_{\mathscr{A}}(g,
C)$ is also epic.}

\vspace{0.2cm}

{\it Proof.} Assume that the given diagram is a pull-back diagram of
$f$ and $g$ and $\Hom_{\mathscr{A}}(C, g)$ is epic. Let $\alpha \in
\Hom_{\mathscr{A}}(C, X)$. Then there exists $\beta \in
\Hom_{\mathscr{A}}(C, N)$ such that $f\alpha=\Hom_{\mathscr{A}}(C,
g)(\beta) =g\beta$. By the universal property of a pull-back
diagram, there exists $\gamma \in \Hom_{\mathscr{A}}(C, M)$ such
that $\alpha=g_1\gamma=\Hom_{\mathscr{A}}(C, g_1)(\gamma)$. So
$\Hom_{\mathscr{A}}(C, g_1)$ is epic and the assertion (1) follows.

Dually, we get the assertion (2). \hfill$\square$

\vspace{0.2cm}

{\bf Lemma 2.5.} {\it (1) Let
$$\xymatrix{K_1 \ar[d]\ar[r]^{h} & K_0\ar[d]^{\alpha_0}\\
W_1 \ar[d]\ar[r]^{g} & W_0\ar[d]\\
X_1 \ar[d]\ar[r]^{f} & X_0\ar[d]\\
0 & 0}$$ be a commutative diagram in
$\mathscr{A}$ with exact columns and $C$ an object in $\mathscr{A}$. If all of
$\Hom_{\mathscr{A}}(h,C)$, $\Hom_{\mathscr{A}}(f,C)$ and
$\Hom_{\mathscr{A}}(\alpha_0,C)$ are epic, then so is
$\Hom_{\mathscr{A}}(g,C)$.

(2) Let
$$\xymatrix{0\ar[d] & 0\ar[d]\\
K_1 \ar[d]\ar[r]^{h} & K_0\ar[d]\\
W_1 \ar[d]^{\beta_1}\ar[r]^{g} & W_0\ar[d]\\
X_1 \ar[r]^{f}& X_0}$$ be a commutative diagram
in $\mathscr{A}$ with exact columns and $C$ an object in $\mathscr{A}$. If all of
$\Hom_{\mathscr{A}}(C,h)$, $\Hom_{\mathscr{A}}(C,f)$ and
$\Hom_{\mathscr{A}}(C,\beta_1)$ are epic, then so is
$\Hom_{\mathscr{A}}(C,g)$.}

\vspace{0.2cm}

{\it Proof.} (1) By assumption, we get the following commutative
diagram with exact columns and rows:
$$\xymatrix{ & & 0\\
0 & \Hom_{\mathscr{A}}(K_1,C)\ar[l] & \Hom_{\mathscr{A}}(K_0,C)\ar[u]\ar[l]_{\Hom_{\mathscr{A}}(h,C)}\\
& \Hom_{\mathscr{A}}(W_1,C)\ar[u] & \Hom_{\mathscr{A}}(W_0,C)
\ar[u]_{\Hom_{\mathscr{A}}(\alpha_0,C)}\ar[l]_{\Hom_{\mathscr{A}}(g,C)}\\
0 & \Hom_{\mathscr{A}}(X_1,C)\ar[u]\ar[l] & \Hom_{\mathscr{A}}(X_0,C)\ar[u]\ar[l]_{\Hom_{\mathscr{A}}(f,C)}\\
& 0\ar[u] & 0\ar[u]}$$ By the snake lemma, $\Hom_{\mathscr{A}}(g,
C)$ is epic.

(2) It is dual to (1). \hfill$\square$

\vspace{0.2cm}

{\bf Definition 2.6.} ([SSW1]) The {\it Gorenstein subcategory}
$\mathcal{G}(\mathscr{C})$ of $\mathscr{A}$ is defined as
$\mathcal{G}(\mathscr{C})=\{M$ is an object in $\mathscr{A}\mid$
there exists an exact sequence: $$\cdots \to C_1 \to C_0 \to C^0 \to
C^1 \to \cdots \eqno{(2.1)}$$ in $\mathscr{C}$, which is both
$\Hom_{\mathscr{A}}(\mathscr{C},-)$-exact and
$\Hom_{\mathscr{A}}(-,\mathscr{C})$-exact, such that $M\cong
\Im(C_0\to C^0)\}$; in this case, $(2.1)$ is called a {\it complete
$\mathscr{C}$-resolution} of $M$.

\vspace{0.2cm}

{\it Remark 2.7.} (1) Let $R$ be a left and right Noetherian ring
and $\mod R$ the category of finitely generated left $R$-modules.
Put $\mathcal{P}(\mod R)$ the subcategory of $\mod
R$ consisting of projective modules. Then
$\mathcal{G}(\mathcal{P}(\mod R))$ coincides with the subcategory of $\mod
R$ consisting of modules with Gorenstein dimension zero ([AB]).

(2) For any ring $R$, $\mathcal{G}(\mathcal{P}(\Mod R))$
(resp. $\mathcal{G}(\mathcal{I}(\Mod R))$)
coincides with the subcategory of $\Mod R$ consisting of Gorenstein
projective (resp. injective) modules ([EJ1]).

(3) Let $R$ be a left Noetherian ring, $S$ a right Noetherian ring
and $_RV_S$ a dualizing bimodule. Put $\mathscr{W}=\{V\bigotimes_SP\mid P$
is projective in $\Mod S\}$ and $\mathscr{U}=\{\Hom_S(V,E)\mid E$
is injective in $\Mod S^{op}\}$. Then $\mathcal{G}(\mathscr{W})$
(resp. $\mathcal{G}(\mathscr{U})$) coincides
with the subcategory of $\Mod R$ consisting of $V$-Gorenstein
projective (resp. injective) modules ([EJL]).

\vspace{0.5cm}

\centerline{\large \bf 3. The constructions of (strongly) proper resolutions}
\centerline{\large \bf  and coproper coresolutions}

\vspace{0.2cm}

In this section, we introduce the notion of strongly (co)proper
(co)resolutions of modules. Then we give a method to construct a
(strongly) proper resolution (resp. coproper coresolution) of the first (resp.
last) term in a short exact sequence from that of the other two
terms, as well as give a method to construct a
(strongly) proper resolution (resp. coproper coresolution) of the last (resp.
first) term in a short exact sequence from that of the other two
terms.

We first give the following easy observation, which is a
generalization of the horseshoe lemma.

\vspace{0.2cm}

{\bf Lemma 3.1.} {\it Let $0\to A \buildrel {f} \over
\longrightarrow A^{'} \buildrel {g} \over \longrightarrow A^{''} \to
0$ be an exact sequence in $\mathscr{A}$.

(1) If there exist morphisms $\alpha\in \Hom_{\mathscr{A}}(C, A)$, $\alpha^{''}\in
\Hom_{\mathscr{A}}(C^{''}, A^{''})$ and $h\in
\Hom_{\mathscr{A}}(C^{''}, A^{'})$ such that $\alpha^{''}=gh$, then
we have the following commutative diagram with exact rows:
$$\xymatrix{0 \ar[r] & C \ar[d]_{\alpha}
\ar[r]^{\binom {1_C} 0} & C\bigoplus C^{''}
\ar@{-->}[d]_{\alpha^{'}}
\ar[r]^{(0, 1_{C^{''}})} & C^{''} \ar[d]_{\alpha^{''}} \ar[r] & 0\\
0 \ar[r] & A \ar[r]^f & A^{'}  \ar[r]^g & A^{''} \ar[r] & 0}$$ where
$\alpha^{'}=(f\alpha, h)(\in \Hom_{\mathscr{A}}(C\bigoplus C^{''}, A^{'}))$.

(2) If there exist morphisms $\beta\in\Hom_{\mathscr{A}}(A, D)$,
$\beta^{''}\in\Hom_{\mathscr{A}}(A^{''}, D^{''})$ and
$k\in\Hom_{\mathscr{A}}(A^{'}, D)$ such that $\beta=kf$, then we
have the following commutative diagram with exact rows:
$$\xymatrix{0 \ar[r] & A \ar[d]_{\beta}\ar[r]^f & A^{'}
\ar@{-->}[d]_{\beta^{'}} \ar[r]^g & A^{''} \ar[d]_{\beta^{''}}\ar[r]& 0\\
0\ar[r] & D \ar[r]^{\binom {1_D} 0} & D\bigoplus D^{''} \ar[r]^{(0,
1_{D^{''}})} & D^{''} \ar[r] & 0}$$ where $\beta^{'}={\binom k
{\beta^{''}g}}(\in \Hom_{\mathscr{A}}(A^{'}, D\bigoplus D^{''}))$.}

\vspace{0.2cm}

The following result contains Theorem 1.2, which provides a method to
construct a (strongly) proper resolution of the first term in a
short exact sequence from that of the last two terms.

\vspace{0.2cm}

{\bf Theorem 3.2.} {\it Let $0\to X \to X^0 \to X^1 \to 0$ be an
exact sequence in $\mathscr{A}$. Let
$$\cdots \to C_i^0 \to \cdots \to C_1^0 \to C_0^0 \to X^0 \to 0 \eqno{(3.1)}$$
be a $\mathscr{C}$-resolution of $X^0$, and let
$$\cdots \to C_i^1 \to \cdots \to C_1^1 \to C_0^1 \to X^1 \to 0 \eqno{(3.2)}$$
be a $\Hom_{\mathscr{A}}(\mathscr{C}, -)$-exact exact sequence in
$\mathscr{A}$. Then

(1) We get the following exact sequences:
$$\cdots \to C_{i+1}^1\bigoplus C_i^0 \to \cdots
\to C_2^1\bigoplus C_1^0 \to C \to X \to 0 \eqno{(3.3)}$$ and
$$0\to C \to C_1^1\bigoplus C_0^0 \to C_0^1 \to 0 \eqno{(3.4)}$$

(2) If both the exact sequences (3.1) and (3.2) are
$\Hom_{\mathscr{A}}(-,\mathscr{C})$-exact, then so is the exact
sequence (3.3).

Assume that $\mathscr{C}$ is closed under finite direct sums and
kernels of epimorphisms. Then we have

(3) If the exact sequence (3.2) is a $\mathscr{C}$-resolution of
$X^1$, then the exact sequence (3.3) is a $\mathscr{C}$-resolution
of $X$.

(4) If both the exact sequences (3.1) and (3.2) are proper
$\mathscr{C}$-resolutions of $X^0$ and $X^1$ respectively, then the
exact sequence (3.3) is a proper $\mathscr{C}$-resolution of $X$.

(5) If both the exact sequences (3.1) and (3.2) are strongly proper
$\mathscr{C}$-resolutions of $X^0$ and $X^1$ respectively, then the
exact sequence (3.3) is a strongly proper $\mathscr{C}$-resolution
of $X$.}

\vspace{0.2cm}

{\it Proof.} (1) Put $K_i^0=\Im(C_i^0\to C_{i-1}^0)$ and
$K_i^1=\Im(C_i^1\to C_{i-1}^1)$ for any $i \geq
1$. Consider the following pull-back diagram:
$$\xymatrix{& & 0 \ar[d] & 0 \ar[d]& &\\
& & K_1^1 \ar@{=}[r] \ar[d] & K_1^1 \ar[d]& &\\
0 \ar[r] & X \ar@{=}[d] \ar[r] & M
\ar[d] \ar[r] &C_0^1 \ar[d] \ar[r] & 0\\
0 \ar[r] & X \ar[r] & X^0 \ar[r] \ar[d] & X^1 \ar[d] \ar[r] & 0 &\\
& & 0 & 0 & & }
$$ Because the third column in the above diagram is $\Hom_{\mathscr{A}}(\mathscr{C}, -)$-exact exact,
so is the middle column by Lemma 2.4(1). Thus by Lemma 3.1(1) we get
the following commutative diagram with exact columns and rows and
the middle row splitting:
$$\xymatrix{& 0 \ar[d] & 0 \ar@{-->}[d] & 0 \ar[d]& &\\
0 \ar@{-->}[r] & K_2^1 \ar[d] \ar@{-->}[r] & W_1 \ar@{-->}[d]
\ar@{-->}[r] & K_1^0 \ar[d] \ar@{-->}[r] & 0\\
0 \ar[r] & C_1^1 \ar[d] \ar[r] & C_1^1\bigoplus C_0^0
\ar@{-->}[d] \ar[r] & C_0^0 \ar[d] \ar[r] & 0\\
0 \ar[r] & K_1^1 \ar[d] \ar[r] & M
\ar@{-->}[d] \ar[r] & X^0 \ar[d] \ar[r] & 0\\
& 0  & 0  & 0 & }
$$ where $W_1=\Ker(C_1^1\bigoplus C_0^0\to M)$. It is easy to verify the upper row in the
above diagram is $\Hom_{\mathscr{A}}(\mathscr{C}, -)$-exact exact.

On the one hand, we have the following pull-back diagram:
$$\xymatrix{& 0 \ar[d] & 0 \ar[d] & &\\
&  W_1 \ar@{=}[r] \ar[d]& W_1 \ar[d] & & \\
0 \ar[r] & C \ar[r] \ar[d] & C_1^1\bigoplus C_0^0 \ar[r] \ar[d] & C_0^1 \ar[r] \ar@{=}[d] & 0\\
0 \ar[r] & X \ar[r] \ar[d] & M \ar[r] \ar[d] & C_0^1 \ar[r] & 0\\
& 0 & 0 & &}$$ On the other hand, again by
Lemma 3.1(1) we get the following commutative diagram with exact
columns and rows and the middle row splitting:

$$\xymatrix{& 0 \ar[d] & 0 \ar@{-->}[d] & 0 \ar[d]& &\\
0 \ar@{-->}[r] & K_3^1 \ar[d] \ar@{-->}[r] & W_2 \ar@{-->}[d]
\ar@{-->}[r] & K_2^0 \ar[d] \ar@{-->}[r] & 0\\
0 \ar[r] & C_2^1 \ar[d] \ar[r] & C_2^1\bigoplus C_1^0
\ar@{-->}[d] \ar[r] & C_1^0 \ar[d] \ar[r] & 0\\
0 \ar[r] & K_2^1 \ar[d] \ar[r] & W_1
\ar@{-->}[d] \ar[r] & K_1^0 \ar[d] \ar[r] & 0\\
& 0  & 0  & 0 & }
$$ where $W_2=\Ker(C_2^1\bigoplus C_1^0\to W_1)$ and the upper row in the
above diagram is $\Hom_{\mathscr{A}}(\mathscr{C}, -)$-exact exact.
Continuing this process, we get the desired exact sequences (3.3)
and (3.4) with $W_i=\Im(C_{i+1}^1\bigoplus C_i^0\to C_i^1\bigoplus
C_{i-1}^0)$ for any $i \geq 2$ and $W_1=\Im(C_2^1\bigoplus C_1^0 \to
C)$.

(2) By assumption, both the first and third columns in the second
diagram are $\Hom_{\mathscr{A}}(-,\mathscr{C})$-exact. Then it is
easy to see that both the first row and the middle column in this
diagram are also $\Hom_{\mathscr{A}}(-,\mathscr{C})$-exact. It
follows that the first column in the third diagram is
$\Hom_{\mathscr{A}}(-,\mathscr{C})$-exact. Also by assumption, the
first and third columns in the last diagram are
$\Hom_{\mathscr{A}}(-,\mathscr{C})$-exact. So the middle column in
this diagram is also $\Hom_{\mathscr{A}}(-,\mathscr{C})$-exact.
Finally, we deduce that the exact sequence (3.3) is
$\Hom_{\mathscr{A}}(-,\mathscr{C})$-exact.

(3) It follows from the assumption and the assertion (1).

(4) Assume that both the exact sequences (3.1) and (3.2) are proper
$\mathscr{C}$-resolutions of $X^0$ and $X^1$ respectively. Then by
the proof of (1) and [EJ2, Lemma 8.2.1], we have that both the
middle column in the second diagram and the first column in the
third diagram are $\Hom_{\mathscr{A}}(\mathscr{C}, -)$-exact exact;
and in particular we have a $\Hom_{\mathscr{A}}(\mathscr{C},
-)$-exact exact sequence:
$$\cdots \to C_{i+1}^1\bigoplus C_i^0 \to \cdots
\to C_2^1\bigoplus C_1^0 \to W_1 \to 0.$$ Thus we get the desired
proper $\mathscr{C}$-resolution of $X$.

(5) If both the exact sequences (3.1) and (3.2) are strongly proper
$\mathscr{C}$-resolutions of $X^0$ and $X^1$ respectively, then
$\Ext_{\mathscr{A}}^1(\mathscr{C}, K_i^j)=0$ for any $i \geq 1$ and
$j=0,1$. By the proof of (1), we have an exact sequence:
$$0\to K_{i+1}^1 \to W_i \to K_i^0 \to 0$$ for any $i \geq 1$.
So $\Ext_{\mathscr{A}}^1(\mathscr{C}, W_i)=0$ for any $i \geq 1$,
and hence the exact sequence (3.3) is a strongly proper
$\mathscr{C}$-resolution of $X$. \hfill$\square$

\vspace{0.2cm}

Based on Theorem 3.2, by using induction on $n$ it is not difficult
to get the following

\vspace{0.2cm}

{\bf Corollary 3.3.} {\it Let $\mathscr{C}$ be closed under finite
direct sums and kernels of epimorphisms, and let $0\to X \to X^0 \to
X^1 \to \cdots \to X^n \to 0$ be an exact sequence in $\mathscr{A}$.
Assume that
$$\cdots \to C_i^j \to \cdots \to C_1^j \to C_0^j \to X^j \to 0 \eqno{(3.5(j))}$$
is a (strongly) proper $\mathscr{C}$-resolution of $X^j$ for any $0
\leq j \leq n$. Then we have

(1) $$\cdots \to \bigoplus _{i=0}^nC_{i+3}^i \to \bigoplus
_{i=0}^nC_{i+2}^i \to \bigoplus _{i=0}^nC_{i+1}^i \to C \to X \to 0
\eqno{(3.6)}$$ is a (strongly) proper $\mathscr{C}$-resolution of
$X$, and there exists an exact sequence:
$$0\to C \to \bigoplus _{i=0}^nC_i^i \to \bigoplus _{i=1}^nC_{i-1}^i
\to \bigoplus _{i=2}^nC_{i-2}^i \to \cdots \to C_0^{n-1}\bigoplus
C_1^n \to C_0^n \to 0.$$

(2) If all $(3.5(j))$ are
$\Hom_{\mathscr{A}}(-,\mathscr{C})$-exact, then so is (3.6).}





\vspace{0.2cm}





The next two results are dual to Theorem 3.2 and Corollary 3.3
respectively. The following result provides a method to construct a
(strongly) coproper coresolution of the last term in a short exact
sequence from that of the first two terms.

\vspace{0.2cm}

{\bf Theorem 3.4.} {\it Let $0\to Y_1 \to Y_0 \to Y \to 0$ be an
exact sequence in $\mathscr{A}$. Let
$$0\to Y_0 \to C_0^0 \to C_0^1 \to \cdots \to C_0^i \to \cdots \eqno{(3.7)}$$
be a $\mathscr{C}$-coresolution of $Y_0$, and let
$$0\to Y_1 \to C_1^0 \to C_1^1 \to \cdots \to C_1^i \to \cdots \eqno{(3.8)}$$
be a $\Hom_{\mathscr{A}}(-, \mathscr{C})$-exact exact sequence in
$\mathscr{A}$.

(1) We get the following exact sequences:
$$0\to Y \to C \to C_0^1\bigoplus C_1^2 \to \cdots
\to C_0^i\bigoplus C_1^{i+1} \to \cdots \eqno{(3.9)}$$ and
$$0\to C_1^0 \to C_0^0\bigoplus C_1^1 \to C \to 0 \eqno{(3.10)}$$

(2) If both the exact sequences (3.7) and (3.8) are
$\Hom_{\mathscr{A}}(\mathscr{C},-)$-exact, then so is the exact
sequence (3.9).

Assume that $\mathscr{C}$ is closed under finite direct sums and
cokernels of monomorphisms. Then we have

(3) If the exact sequence (3.8) is a $\mathscr{C}$-coresolution of
$Y_1$, then the exact sequence (3.9) is a $\mathscr{C}$-coresolution
of $Y$.

(4) If both the exact sequences (3.7) and (3.8) are coproper
$\mathscr{C}$-coresolutions of $Y_0$ and $Y_1$ respectively, then
the exact sequence (3.9) is a coproper $\mathscr{C}$-coresolution of
$Y$.

(5) If both the exact sequences (3.7) and (3.8) are strongly
coproper $\mathscr{C}$-coresolutions of $Y_0$ and $Y_1$
respectively, then the exact sequence (3.9) is a strongly coproper
$\mathscr{C}$-coresolution of $Y$.}

\vspace{0.2cm}

{\it Proof.} It is dual to the proof of Theorem 3.2, we give the
proof here for the sake of completeness.

(1) Put $K_0^i=\Im(C_0^{i-1}\to C^i_0)$ and $K_1^i=\Im(C^{i-1}_1\to
C^i_1)$ for any $i \geq 1$. Consider the following push-out diagram:
$$\xymatrix{& 0 \ar[d] & 0 \ar[d]& & &\\
0 \ar[r] & Y_1 \ar[r] \ar[d] & Y_0 \ar[d] \ar[r] & Y
\ar@{=}[d] \ar[r] & 0\\
0 \ar[r] & C^0_1 \ar[r] \ar[d] & N \ar[r] \ar[d] & Y \ar[r] & 0\\
& K_1^1 \ar[d] \ar@{=}[r] & K_1^1 \ar[d] & & &\\
& 0 & 0 & & & }
$$
Because the first column in the above diagram is
$\Hom_{\mathscr{A}}(-, \mathscr{C})$-exact exact, so is the middle
column by Lemma 2.4(2). Then by Lemma 3.1(2) we get the following
commutative diagram with exact columns and rows and the middle row
splitting:
$$\xymatrix{& 0 \ar[d] & 0 \ar@{-->}[d] & 0 \ar[d]& &\\
0 \ar[r] & Y_0 \ar[d] \ar[r] & N
\ar@{-->}[d] \ar[r] & K_1^1 \ar[d] \ar[r] & 0\\
0 \ar[r] & C_0^0 \ar[d] \ar[r] & C_0^0\bigoplus C_1^1
\ar@{-->}[d] \ar[r] & C_1^1 \ar[d] \ar[r] & 0\\
0 \ar@{-->}[r] & K_0^1 \ar[d] \ar@{-->}[r] & W^1 \ar@{-->}[d]
\ar@{-->}[r] & K_1^2 \ar[d] \ar@{-->}[r] & 0\\
& 0  & 0  & 0 & }
$$ where $W^1=\Coker(N \to C_0^0\bigoplus C_1^1)$. It is easy to verify that
the bottom row in the above diagram is $\Hom_{\mathscr{A}}(-,
\mathscr{C})$-exact exact.

On the one hand, we have the following push-out diagram:
$$\xymatrix{& & 0 \ar[d] & 0 \ar[d] & \\
0 \ar[r] & C_1^0 \ar[r] \ar@{=}[d]
& N \ar[r] \ar[d] & Y \ar[r] \ar[d] & 0\\
0 \ar[r] & C_1^0 \ar[r] & C_0^0\bigoplus C_1^1 \ar[r] \ar[d] & C \ar[r] \ar[d] & 0\\
& & W^1 \ar@{=}[r] \ar[d] & W^1 \ar[d] &\\
& & 0 & 0 &}$$ On the other hand, again by Lemma
3.1(2) we get the following commutative diagram with exact columns
and rows and the middle row splitting:
$$\xymatrix{& 0 \ar[d] & 0 \ar@{-->}[d] & 0 \ar[d]& &\\
0 \ar[r] & K_0^1 \ar[d] \ar[r] & W^1
\ar@{-->}[d] \ar[r] & K_1^2 \ar[d] \ar[r] & 0\\
0 \ar[r] & C_0^1 \ar[d] \ar[r] & C_0^1\bigoplus C_1^2
\ar@{-->}[d] \ar[r] & C_1^2 \ar[d] \ar[r] & 0\\
0 \ar@{-->}[r] & K_0^2 \ar[d] \ar@{-->}[r] & W^2 \ar@{-->}[d]
\ar@{-->}[r] & K_1^3 \ar[d] \ar@{-->}[r] & 0\\
& 0  & 0  & 0 & }
$$
where $W^2=\Coker(W^1 \to C_0^1\bigoplus C_1^2)$ and the bottom row
in the above diagram is $\Hom_{\mathscr{A}}(-, \mathscr{C})$-exact
exact. Continuing this process, we get the desired exact sequences
(3.9) and (3.10) with $W^i=\Im(C^{i-1}_0\bigoplus C^i_1\to
C^i_0\bigoplus C^{i+1}_1)$ for any $i \geq 2$ and $W^1=\Im(C\to
C_0^1\bigoplus C_1^2)$.

(2) By assumption, both the first and third columns in the second
diagram are $\Hom_{\mathscr{A}}(\mathscr{C},-)$-exact. Then it is
easy to see that both the first row and the middle column in this
diagram are also $\Hom_{\mathscr{A}}(\mathscr{C},-)$-exact. It
follows that the first column in the third diagram is
$\Hom_{\mathscr{A}}(\mathscr{C},-)$-exact. Also by assumption, the
first and third columns in the last diagram are
$\Hom_{\mathscr{A}}(\mathscr{C},-)$-exact. So the middle column in
this diagram is also $\Hom_{\mathscr{A}}(\mathscr{C},-)$-exact.
Finally, we deduce that the exact sequence (3.9) is
$\Hom_{\mathscr{A}}(\mathscr{C},-)$-exact.

(3) It follows from the assumption and the assertion (1).

(4) Assume that both the exact sequences (3.7) and (3.8) are
coproper $\mathscr{C}$-coresolutions of $X_0$ and $X_1$
respectively. Then by the proof of (1) and the dual version of [EJ2,
Lemma 8.2.1], we have that both the middle column in the second
diagram and the first column in the third diagram are
$\Hom_{\mathscr{A}}(-, \mathscr{C})$-exact exact; and in particular
we have a $\Hom_{\mathscr{A}}(-, \mathscr{C})$-exact exact sequence:
$$0\to W^1 \to C_0^1\bigoplus C_1^2 \to \cdots
\to C_0^i\bigoplus C_1^{i+1} \to \cdots .$$ Thus we get the desired
coproper $\mathscr{C}$-coresolution of $X$.

(5) If both the exact sequences (3.7) and (3.8) are strongly
coproper $\mathscr{C}$-coresolutions of $Y_0$ and $Y_1$
respectively, then $\Ext_{\mathscr{A}}^1(K^i_j, \mathscr{C})=0$ for
any $i \geq 1$ and $j=0,1$. By the proof of (1), we have an exact
sequence:
$$0\to K_0^i \to W^i \to K_1^{i+1} \to 0$$ for any $i \geq 1$.
So $\Ext_{\mathscr{A}}^1(W^i, \mathscr{C})=0$ for any $i \geq 1$,
and hence the exact sequence (3.9) is a strongly coproper
$\mathscr{C}$-coresolution of $Y$. \hfill$\square$

\vspace{0.2cm}

Based on Theorem 3.4, by using induction on $n$ it is not difficult
to get the following

\vspace{0.2cm}

{\bf Corollary 3.5.} {\it Let $\mathscr{C}$ be closed under finite
direct sums and cokernels of monomorphisms and let $0\to Y_n \to
\cdots \to Y_1 \to Y_0 \to Y \to 0$ be an exact sequence in
$\mathscr{A}$. Assume that
$$0\to Y_j \to C_j^0 \to C_j^1 \to \cdots \to C_j^i \to \cdots \eqno{(3.11(j))}$$
is a (strongly) coproper $\mathscr{C}$-coresolution of $Y_j$ for any
$0 \leq j \leq n$. Then we have

(1) $$0\to Y \to C \to \bigoplus _{i=0}^nC_i^{i+1} \to \bigoplus
_{i=0}^nC_i^{i+2} \to \bigoplus _{i=0}^nC_i^{i+3} \to \cdots
\eqno{(3.12)}$$ is a (strongly) coproper $\mathscr{C}$-coresolution
of $Y$, and there exists an exact sequence:
$$0\to C_n^0 \to C_{n-1}^0\bigoplus
C_n^1 \to \cdots \to \bigoplus _{i=2}^nC_i^{i-2} \to \bigoplus
_{i=1}^nC_i^{i-1} \to \bigoplus _{i=0}^nC_i^i \to C \to 0.$$

(2) If all $(3.11(j))$ are
$\Hom_{\mathscr{A}}(\mathscr{C},-)$-exact, then so is (3.12).}

\vspace{0.2cm}

The following result contains Theorem 1.3, which provides a method to
construct a (strongly) proper resolution of the last term in a short
exact sequence from that of the first two terms.

\vspace{0.2cm}

{\bf Theorem 3.6.} {\it Let
$$0\to X_1 \to X_0 \to X \to 0 \eqno{(3.13)}$$ be an exact sequence
in $\mathscr{A}$. Let
$$C_0^n \to \cdots \to C_0^1 \to C_0^0 \to X_0 \to 0 \eqno{(3.14)}$$
be a $\Hom_{\mathscr{A}}(\mathscr{C},-)$-exact exact sequence and
$$C_1^{n-1} \to \cdots \to C_1^1 \to C_1^0 \to X_1 \to 0 \eqno{(3.15)}$$
a $\mathscr{C}$-resolution of $X_1$ in $\mathscr{A}$.

(1) We get the following exact sequence:
$$C_0^n\bigoplus C_1^{n-1} \to \cdots
\to C_0^2\bigoplus C_1^1 \to C_0^1\bigoplus C_1^0 \to C_0^0 \to X
\to 0 \eqno{(3.16)}$$

(2) If all the exact sequences (3.13)--(3.15) are
$\Hom_{\mathscr{A}}(-,\mathscr{C})$-exact, then so is the exact
sequence (3.16).

Assume that $\mathscr{C}$ is closed under finite direct sums. Then
we have

(3) If the exact sequence (3.14) is a $\mathscr{C}$-resolution of
$X_0$, then the exact sequence (3.16) is a $\mathscr{C}$-resolution
of $X$.

(4) If the exact sequence (3.13) is
$\Hom_{\mathscr{A}}(\mathscr{C},-)$-exact and both the exact
sequences (3.14) and (3.15) are proper $\mathscr{C}$-resolutions of
$X_0$ and $X_1$ respectively, then the exact sequence (3.16) is a
proper $\mathscr{C}$-resolution of $X$.

(5) If the exact sequence (3.13) is strongly
$\Hom_{\mathscr{A}}(\mathscr{C},-)$-exact and both the exact
sequences (3.14) and (3.15) are strongly proper
$\mathscr{C}$-resolutions of $X_0$ and $X_1$ respectively, then the
exact sequence (3.16) is a strongly proper $\mathscr{C}$-resolution
of $X$.}

\vspace{0.2cm}

{\it Proof.} (1) Put $K_j^i=\Im(C_j^i\to C_j^{i-1})$ for any $1\leq i \leq
n-j$ and $j=0,1$. Consider the following pull-back diagram:
$$\xymatrix{& 0 \ar[d] & 0 \ar[d] & & \\
& K_0^1\ar[d]\ar@{=}[r] & K_0^1\ar[d] & & \\
0 \ar[r] & W_1 \ar[r]\ar[d] & C_0^0 \ar[r]\ar[d]
& X \ar[r]\ar@{=}[d] & 0\\
0 \ar[r] & X_1 \ar[r]\ar[d] & X_0 \ar[r]\ar[d]
& X \ar[r] & 0\\
& 0 & 0 & & }$$ Note that the middle column in the above diagram is
$\Hom_{\mathscr{A}}(\mathscr{C},-)$-exact exact. So by Lemma 2.4(1),
the first column is also $\Hom_{\mathscr{A}}(\mathscr{C},-)$-exact
exact. Then by Lemma 3.1(1) we get the following commutative diagram
with exact columns and rows and the middle row splitting:

$$\xymatrix{& 0 \ar[d] & 0 \ar@{-->}[d] & 0 \ar[d]& &\\
0 \ar@{-->}[r] & K_0^2 \ar[d] \ar@{-->}[r] & W_2 \ar@{-->}[d]
\ar@{-->}[r] & K_1^1 \ar[d] \ar@{-->}[r] & 0\\
0 \ar[r] & C_0^1 \ar[d] \ar[r] & C_0^1\bigoplus C_1^0
\ar@{-->}[d] \ar[r] & C_1^0 \ar[d] \ar[r] & 0\\
0 \ar[r] & K_0^1 \ar[d] \ar[r] & W_1
\ar@{-->}[d] \ar[r] & X_1 \ar[d] \ar[r] & 0\\
& 0  & 0  & 0 & }
$$ where $W_2=\Ker(C_0^1\bigoplus C_1^0\to W_1)$.
It is easy to check that the upper row in the above diagram is
$\Hom_{\mathscr{A}}(\mathscr{C},-)$-exact exact. Then by using Lemma
3.1(1) iteratively we get the exact sequence (3.16) with
$W_i=\Im(C_0^i\bigoplus C_1^{i-1}\to C_0^{i-1}\bigoplus C_1^{i-2})$
for any $2\leq i \leq n$ and $W_1=\Im(C_0^1\bigoplus C_1^0 \to
C_0^0)$.

(2) By assumption, both the third row and the middle column in the first diagram
are $\Hom_{\mathscr{A}}(-,\mathscr{C})$-exact. Then the middle row in this
diagram is $\Hom_{\mathscr{A}}(-,\mathscr{C})$-exact by Lemma
2.5(1). It is easy to see that the first column in this diagram is
also $\Hom_{\mathscr{A}}(-,\mathscr{C})$-exact.
Both the first and third columns in the second diagram are $\Hom_{\mathscr{A}}(-,\mathscr{C})$-exact
also by assumption, so the middle column in this diagram is
also $\Hom_{\mathscr{A}}(-,\mathscr{C})$-exact. Finally, we deduce that
the exact sequence (3.16) is
$\Hom_{\mathscr{A}}(-,\mathscr{C})$-exact.

(3) It follows from the assumption and the assertion (1).

(4) If the exact sequence (3.13) is
$\Hom_{\mathscr{A}}(\mathscr{C},-)$-exact and both the exact
sequences (3.14) and (3.15) are proper $\mathscr{C}$-resolutions of
$X_0$ and $X_1$ respectively, then the middle row in the first
diagram is $\Hom_\mathscr{A}{}(\mathscr{C},-)$-exact exact by Lemma
2.4(1). Thus by [EJ2, Lemma 8.2.1], the exact sequence (3.16) is a
proper $\mathscr{C}$-resolution of $X$.

(5) If the exact sequence (3.13) is strongly
$\Hom_{\mathscr{A}}(\mathscr{C},-)$-exact and both the exact
sequences (3.14) and (3.15) are strongly proper
$\mathscr{C}$-resolutions of $X_0$ and $X_1$ respectively, then
$\Ext_{\mathscr{A}}^1(\mathscr{C}, K_j^i)=0$ for any $1\leq i \leq
n-j$ and $j=0,1$. By the proof of (1), we have an exact sequence:
$$0\to K_0^i \to W_i \to K_1^{i-1} \to 0$$ for any $1\leq i \leq n$ (where $K_1^0=X_1$).
So $\Ext_{\mathscr{A}}^1(\mathscr{C}, W_i)=0$ for any $1\leq i \leq
n$, and hence the exact sequence (3.16) is a strongly proper
$\mathscr{C}$-resolution of $X$. \hfill$\square$

\vspace{0.2cm}

Based on Theorem 3.6, by using induction on $n$ it is not difficult
to get the following

\vspace{0.2cm}

{\bf Corollary 3.7.} {\it Let $\mathscr{C}$ be closed under finite
direct sums, and let $$X_n \to \cdots \to X_1 \to X_0 \to X \to 0
\eqno{(3.17)}$$ and
$$C_j^{n-j} \to \cdots \to C_j^1 \to C_j^0 \to X_j \to 0 \eqno{(3.18(j))}$$
be exact sequences in $\mathscr{A}$ for any $0 \leq j \leq n$.

(1) Let the exact sequence (3.17) be
$\Hom_{\mathscr{A}}(\mathscr{C},-)$-exact and $(3.18(j))$ a
proper $\mathscr{C}$-resolution of $X_j$ for any $0 \leq j \leq n$.
Then $$\bigoplus _{i=0}^nC_i^{n-i}\to \bigoplus
_{i=0}^{n-1}C_i^{(n-1)-i} \to \cdots \to C_0^1\bigoplus C_1^0 \to
C_0^0 \to X \to 0 \eqno{(3.19)}$$ is a proper
$\mathscr{C}$-resolution of $X$; furthermore, if (3.17) and all $(3.18(j))$ are
$\Hom_{\mathscr{A}}(-,\mathscr{C})$-exact, then so is (3.19).

(2) Let the exact sequence (3.17) be strongly
$\Hom_{\mathscr{A}}(\mathscr{C},-)$-exact and $(3.18(j))$ a
strongly proper $\mathscr{C}$-resolution of $X_j$ for any $0\leq j
\leq n$. Then (3.19) is a strongly proper $\mathscr{C}$-resolution
of $X$.}





\vspace{0.2cm}

The next two results are dual to Theorem 3.6 and Corollary 3.7
respectively. The following result provides a method to construct a
(strongly) coproper coresolution of the first term in a short exact
sequence from that of the last two terms.

\vspace{0.2cm}

{\bf Theorem 3.8.} {\it Let
$$0\to Y \to Y^0 \to Y^1 \to 0 \eqno{(3.20)}$$ be an exact sequence
in $\mathscr{A}$. Let
$$0\to Y^0 \to C_0^0 \to C_1^0 \to \cdots \to C_n^0 \eqno{(3.21)}$$
be a $\Hom_{\mathscr{A}}(-, \mathscr{C})$-exact exact sequence and
$$0\to Y^1 \to C_0^1 \to C_1^1 \to \cdots \to C_{n-1}^1 \eqno{(3.22)}$$
a $\mathscr{C}$-coresolution of $Y^1$ in $\mathscr{A}$.

(1) We get the following exact sequence:
$$0\to Y \to C_0^0 \to C_0^1\bigoplus C_1^0 \to  C_1^1\bigoplus C_2^0 \to \cdots
\to C_{n-1}^1\bigoplus C_n^0 \eqno{(3.23)}$$

(2) If all the exact sequences (3.20)--(3.22) are
$\Hom_{\mathscr{A}}(\mathscr{C},-)$-exact, then so is the exact
sequence (3.23).

Assume that $\mathscr{C}$ is closed under finite direct sums. Then
we have

(3) If the exact sequence (3.21) is a $\mathscr{C}$-coresolution of
$Y^0$, then the exact sequence (3.23) is a
$\mathscr{C}$-coresolution of $X$.

(4) If the exact sequence (3.20) is $\Hom_{\mathscr{A}}(-,
\mathscr{C})$-exact and both the exact sequences (3.21) and (3.22)
are coproper $\mathscr{C}$-coresolutions of $Y^0$ and $Y^1$
respectively, then the exact sequence (3.23) is a coproper
$\mathscr{C}$-coresolution of $Y$.

(5) If the exact sequence (3.20) is strongly $\Hom_{\mathscr{A}}(-,
\mathscr{C})$-exact and both the exact sequences (3.21) and (3.22)
are strongly coproper $\mathscr{C}$-coresolutions of $Y^0$ and $Y^1$
respectively, then the exact sequence (3.23) is a strongly coproper
$\mathscr{C}$-coresolution of $Y$.}

\vspace{0.2cm}

{\it Proof.} It is dual to the proof of Theorem 3.6, we give the
proof here for the sake of completeness.

(1) Put $K_i^j=\Im(C_{i-1}^j\to C_i^j)$ for any $1\leq i\leq n-j$
and $j=0,1$. Consider the following push-out diagram:
$$\xymatrix{& & 0 \ar[d] & 0 \ar[d] & \\
0 \ar[r] & Y \ar[r] \ar@{=}[d]
& Y^0 \ar[r] \ar[d] & Y^1 \ar[r] \ar[d] & 0\\
0 \ar[r] & Y \ar[r] & C_0^0 \ar[r] \ar[d] & W^1 \ar[r] \ar[d] & 0\\
& & K_1^0 \ar@{=}[r] \ar[d] & K_1^0 \ar[d] &\\
& & 0 & 0 &}$$ Note that the middle column in the above diagram is
$\Hom_{\mathscr{A}}(-, \mathscr{C})$-exact exact by assumption. So
the third column is also $\Hom_{\mathscr{A}}(-, \mathscr{C})$-exact
exact by Lemma 2.4(2). Then by Lemma 3.1(2) we get the following
commutative diagram with exact columns and rows and the middle row
splitting:
$$\xymatrix{& 0 \ar[d] & 0 \ar@{-->}[d] & 0 \ar[d]& &\\
0 \ar[r] & Y^1 \ar[d] \ar[r] & W^1
\ar@{-->}[d] \ar[r] & K_1^0 \ar[d] \ar[r] & 0\\
0 \ar[r] & C_0^1 \ar[d] \ar[r] & C_0^1\bigoplus C_1^0
\ar@{-->}[d] \ar[r] & C_1^0 \ar[d] \ar[r] & 0\\
0 \ar@{-->}[r] & K_1^1 \ar[d] \ar@{-->}[r] & W^2 \ar@{-->}[d]
\ar@{-->}[r] & K_2^0 \ar[d] \ar@{-->}[r] & 0\\
& 0  & 0  & 0 & }$$ where $W^2=\Coker(W^1 \to C_0^1\bigoplus
C_1^0)$. It is easy to verify that the bottom row in the above
diagram is $\Hom_{\mathscr{A}}(-, \mathscr{C})$-exact exact. Then by
using Lemma 3.1(2) iteratively we get the exact sequence (3.23) with
$W^i=\Im(C_{i-2}^1\bigoplus C_{i-1}^0 \to C_{i-1}^1\bigoplus C_i^0)$
for any $2\leq i \leq n$ and $W^1=\Im(C_0^0 \to C_0^1\bigoplus
C_1^0)$.

(2) By assumption, both the first row and the middle column in the first diagram are
$\Hom_{\mathscr{A}}(\mathscr{C},-)$-exact. Then the middle row in this
diagram is $\Hom_{\mathscr{A}}(\mathscr{C},-)$-exact by Lemma
2.5(2). It is easy to see that the third column in this diagram is
also $\Hom_{\mathscr{A}}(\mathscr{C},-)$-exact.
Both the first and third columns in the second diagram are $\Hom_{\mathscr{A}}(\mathscr{C},-)$-exact
also by assumption, so the middle column in this diagram is
also $\Hom_{\mathscr{A}}(\mathscr{C},-)$-exact. Finally, we deduce that
the exact sequence (3.23) is $\Hom_{\mathscr{A}}(\mathscr{C},-)$-exact.

(3) It follows from the assumption and the assertion (1).

(4) If the exact sequence (3.20) is $\Hom_{\mathscr{A}}(-,
\mathscr{C})$-exact and both the exact sequences (3.21) and (3.22)
are coproper $\mathscr{C}$-coresolutions of $Y^0$ and $Y^1$
respectively, then the middle row in the first diagram is
$\Hom_{\mathscr{A}}(-, \mathscr{C})$-exact exact by Lemma 2.4(2).
Thus by the dual version of [EJ2, Lemma 8.2.1], the exact sequence
(3.23) is a coproper $\mathscr{C}$-coresolution of $Y$.

(5) If the exact sequence (3.20) is strongly $\Hom_{\mathscr{A}}(-,
\mathscr{C})$-exact and both the exact sequences (3.21) and (3.22)
are strongly coproper $\mathscr{C}$-coresolutions of $Y^0$ and $Y^1$
respectively, then $\Ext_{\mathscr{A}}^1(K_i^j, \mathscr{C})=0$ for
any $1\leq i\leq n-j$ and $j=0,1$. By the proof of (1), we have an
exact sequence:
$$0\to K_{i-1}^1 \to W^i \to K_i^0 \to 0$$ for any $1\leq i \leq n$ (where $K_0^1=Y^1$).
So $\Ext_{\mathscr{A}}^1(W^i, \mathscr{C})=0$ for any $1\leq i \leq
n$, and hence the exact sequence (3.23) is a strongly coproper
$\mathscr{C}$-coresolution of $Y$. \hfill$\square$

\vspace{0.2cm}

Based on Theorem 3.8, by using induction on $n$ it is not difficult
to get the following

\vspace{0.2cm}

{\bf Corollary 3.9.} {\it Let $\mathscr{C}$ be closed under finite
direct sums, and let $$0\to Y \to Y^0 \to Y^1 \to \cdots \to Y^n
\eqno{(3.24)}$$ and
$$0\to Y^j \to C^j_0 \to C^j_1 \to \cdots \to C^j_{n-j} \eqno{(3.25(j))}$$
be exact sequences in $\mathscr{A}$ for any $0 \leq j \leq n$.

(1) Let the exact sequence (3.24) be
$\Hom_{\mathscr{A}}(-,\mathscr{C})$-exact and $(3.25(j))$ a
coproper $\mathscr{C}$-coresolution of $Y^j$ for any $0\leq j \leq
n$. Then
$$0\to Y \to C_0^0 \to C_1^0\bigoplus C_0^1 \to \cdots \to
\bigoplus_{i=0}^{n-1}C^i_{(n-1)-i} \to \bigoplus
_{i=0}^nC^i_{n-i}\eqno{(3.26)}$$ is a coproper
$\mathscr{C}$-coresolution of $Y$; furthermore, if (3.24) and all $(3.25(j))$
are $\Hom_{\mathscr{A}}(\mathscr{C},-)$-exact, then so is (3.26).

(2) Let the exact sequence (3.24) be strongly
$\Hom_{\mathscr{A}}(-,\mathscr{C})$-exact and $(3.25(j))$ a
strongly coproper $\mathscr{C}$-coresolution of $Y^j$ for any $0\leq
j \leq n$. Then (3.26) is a strongly coproper
$\mathscr{C}$-coresolution of $Y$.}





\vspace{0.5cm}

\centerline{\large \bf 4. Gorenstein categories}

\vspace{0.2cm}

In the rest of this paper, all subcategories are additive
subcategories of $\mathscr{A}$, that is, all subcategories are
closed under finite direct sums. Set
$\mathcal{G}^0(\mathscr{C})=\mathscr{C}$,
$\mathcal{G}^1(\mathscr{C})=\mathcal{G}(\mathscr{C})$, and
inductively set
$\mathcal{G}^{n+1}(\mathscr{C})=\mathcal{G}(\mathcal{G}^n(\mathscr{C}))$
for any $n\geq 1$. As an application of the results in the above
section, we get the the following result, which answers Question 1.1
affirmatively.

\vspace{0.2cm}

{\bf Theorem 4.1.} {\it
$\mathcal{G}^n(\mathscr{C})=\mathcal{G}(\mathscr{C})$ for any $n\geq
1$.}

\vspace{0.2cm}

{\it Proof.} It is easy to see that
$\mathscr{C}\subseteq\mathcal{G}(\mathscr{C})\subseteq
\mathcal{G}^2(\mathscr{C})\subseteq\mathcal{G}^3(\mathscr{C})\subseteq\cdots$
is an ascending chain of additive subcategories of $\mathscr{A}$.

Let $M$ be an object in $\mathcal{G}^2(\mathscr{C})$ and
$$\cdots\to G_1 \to G_0 \to G^0 \to G^1 \to \cdots$$ a complete
$\mathcal{G}(\mathscr{C})$-resolution of $M$ with $M\cong \Im(G_0 \to G^0)$.
Then for any $j\geq 0$, there exist exact sequences:
$$\cdots \to C_j^{i} \to \cdots \to C_j^1 \to C_j^0 \to G_j \to 0$$
and
$$0\to G^j \to B^j_0 \to B^j_1 \to \cdots \to B^j_{i} \to \cdots$$
in $\mathscr{A}$ with all $C_j^{i}$ and $B^j_{i}$ objects in
$\mathscr{C}$, which
are both $\Hom_{\mathscr{A}}(\mathscr{C},-)$-exact and
$\Hom_{\mathscr{A}}(-,\mathscr{C})$-exact. By Corollaries 3.7 and
3.9, we get exact sequences:
$$\cdots \to \bigoplus _{i=0}^nC_i^{n-i}\to \cdots \to C_0^1\bigoplus C_1^0 \to C_0^0 \to M \to 0$$
and
$$0\to M \to B_0^0 \to B_1^0\bigoplus B_0^1 \to \cdots \to \bigoplus_{i=0}^nB^i_{n-i}\to \cdots$$
which are both $\Hom_{\mathscr{A}}(\mathscr{C},-)$-exact and
$\Hom_{\mathscr{A}}(-,\mathscr{C})$-exact. So
$$\cdots \to \bigoplus _{i=0}^nC_i^{n-i}\to \cdots \to C_0^1\bigoplus C_1^0 \to C_0^0 \to
B_0^0 \to B_1^0\bigoplus B_0^1 \to \cdots \to
\bigoplus_{i=0}^nB^i_{n-i}\to \cdots$$ is a complete
$\mathscr{C}$-resolution of $M$ with $M\cong \Im(C^0_0 \to B^0_0)$, and hence $M$ is an object in
$\mathcal{G}(\mathscr{C})$ and $\mathcal{G}^2(\mathscr{C})\subseteq
\mathcal{G}(\mathscr{C})$. Thus we have that
$\mathcal{G}^2(\mathscr{C})=\mathcal{G}(\mathscr{C})$. By using
induction on $n$ we get easily the assertion. \hfill$\square$

\vspace{0.2cm}

{\bf Definition 4.2.} (cf. [SSW1]) A subcategory $\mathscr{X}$ of
$\mathscr{C}$ is called a {\it generator} (resp. {\it cogenerator})
for $\mathscr{C}$ if for any object $C$ in $\mathscr{C}$, there
exists an exact sequence $0\to C^{'} \to X \to C \to 0$ (resp. $0\to
C \to X \to C^{'} \to 0$) in $\mathscr{C}$ with $X$ an object in
$\mathscr{X}$; and $\mathscr{X}$ is called a {\it projective
generator} (resp. an {\it injective cogenerator}) for $\mathscr{C}$
if $\mathscr{X}$ is a generator (resp. cogenerator) for
$\mathscr{C}$ and $\Ext_{\mathscr{A}}^i(X,C)=0$ (resp.
$\Ext_{\mathscr{A}}^i(C,X)=0)$ for any object $X$ in $\mathscr{X}$,
any object $C$ in $\mathscr{C}$ and $i\geq 1$.

\vspace{0.2cm}

As an immediate consequence of Theorem 4.1, we get the following
three corollaries. The first one generalizes [SSW1, Proposition 4.6]
and answers positively a question in [SSW1, p.492].

\vspace{0.2cm}

{\bf Corollary 4.3.} {\it (1) If $\mathscr{X}$ is a (projective)
generator for $\mathscr{C}$, then $\mathscr{X}$ is a (projective)
generator for $\mathcal{G}^n(\mathscr{C})$.

(2) If $\mathscr{X}$ is an (injective) cogenerator for
$\mathscr{C}$, then $\mathscr{X}$ is an (injective) cogenerator for
$\mathcal{G}^n(\mathscr{C})$.}

\vspace{0.2cm}

We define $\res\widetilde{\mathscr{C}}=\{M$ is an object in
$\mathscr{A}\mid$ there exists a
$\Hom_{\mathscr{A}}(\mathscr{C},-)$-exact exact sequence $\cdots \to C_i
\to \cdots \to C_1 \to C_0 \to M\to 0$ in $\mathscr{A}$ with all
$C_i$ objects in $\mathscr{C}\}$. Dually, we define
$\cores\widetilde{\mathscr{C}}=\{M$ is an object in
$\mathscr{A}\mid$ there exists a
$\Hom_{\mathscr{A}}(-,\mathscr{C})$-exact exact sequence $0 \to C^0 \to
C^1\to \cdots \to C^i \to \cdots$ in $\mathscr{A}$ with all $C^i$
objects in $\mathscr{C}\}$ (see [SSW1]).

The next corollary generalizes [SSW1, Theorem 4.8].

\vspace{0.2cm}

{\bf Corollary 4.4.} {\it Let $\mathscr{X}$ be a subcategory of
$\mathscr{C}$. Then we have

(1) If $\mathscr{C}\subseteq\res\widetilde{\mathscr{X}}$, then
$\mathcal{G}^n(\mathscr{C})\subseteq\res\widetilde{\mathscr{X}}$ for
any $n\geq 0$.

(2) If $\mathscr{C}\subseteq\cores\widetilde{\mathscr{X}}$, then
$\mathcal{G}^n(\mathscr{C})\subseteq\cores\widetilde{\mathscr{X}}$
for any $n\geq 0$.}

\vspace{0.2cm}

{\it Proof.} (1) Let
$\mathscr{C}\subseteq\res\widetilde{\mathscr{X}}$. Because
$\mathscr{X}$ is closed under finite sums,
$\mathcal{G}(\mathscr{C})\subseteq\res\widetilde{\mathscr{X}}$ by
Corollary 3.7. Then the assertion follows from Theorem 4.1.

(2) It is dual to (1). \hfill$\square$

\vspace{0.2cm}

We also have the following corollary, which generalizes [SSW1,
Theorem 4.9].

\vspace{0.2cm}

{\bf Corollary 4.5.} {\it Let $\mathscr{X}$ be a subcategory of
$\mathscr{C}$. If
$\mathscr{C}\subseteq\res\widetilde{\mathscr{X}}\bigcap
\cores\widetilde{\mathscr{X}}$, then
$\mathcal{G}^n(\mathscr{C})\subseteq\mathcal{G}(\mathscr{X})$ for
any $n\geq 1$.}

\vspace{0.2cm}

{\it Proof.} By using an argument similar to that in the proof of
Theorem 4.1, we get
$\mathcal{G}(\mathscr{C})\subseteq\mathcal{G}(\mathscr{X})$. Then
the assertion follows from Theorem 4.1. \hfill$\square$

\vspace{0.2cm}

As another application of the results in the above section, we get
the the following result, in which the second assertion shows that
the assumption of the self-orthogonality of $\mathscr{C}$ in [SSW1,
Proposition 4.11] is superfluous.

\vspace{0.2cm}

{\bf Theorem 4.6.} {\it (1) Both $\res\widetilde{\mathscr{C}}$ and
$\cores\widetilde{\mathscr{C}}$ are closed under direct summands.

(2) $\mathcal{G}(\mathscr{C})$ is closed under direct summands.}

\vspace{0.2cm}

{\it Proof.} Assume that $M=X\bigoplus Y$ and
$$0\to Y \buildrel {{\binom 0
{1_Y}}} \over \longrightarrow M \buildrel {(1_X,0)} \over
\longrightarrow X \to 0$$ is an exact and split sequence.

(1) We only prove $\res\widetilde{\mathscr{C}}$ is closed under
direct summands. Dually, we get that $\cores\widetilde{\mathscr{C}}$
is also closed under direct summands.

Let $M$ be an object in $\res\widetilde{\mathscr{C}}$ and
$$\cdots \buildrel {f_{i+1}} \over \longrightarrow C_i \buildrel {f_{i}} \over \longrightarrow
\cdots \buildrel {f_{3}} \over \longrightarrow C_2 \buildrel {f_{2}}
\over \longrightarrow C_1 \buildrel {f_{1}} \over \longrightarrow
C_0 \buildrel {f_{0}} \over \longrightarrow M\to 0$$ a
$\Hom_{\mathscr{A}}(\mathscr{C},-)$-exact exact sequence in
$\mathscr{A}$ with all $C^i$ objects in $\mathscr{C}$. Then
$$C_0 \buildrel {(1_X,0) f_{0}} \over \longrightarrow X \to 0$$ is a
$\Hom_{\mathscr{A}}(\mathscr{C},-)$-exact exact sequence. Similarly,
$$C_0 \buildrel {(1_Y,0) f_{0}} \over \longrightarrow Y \to 0$$ is a
$\Hom_{\mathscr{A}}(\mathscr{C},-)$-exact exact sequence. By Theorem
3.6, we get the following $\Hom_{\mathscr{A}}(\mathscr{C},-)$-exact
exact sequences:
$$C_0\bigoplus C_1 \to C_0 \to X \to 0$$
and
$$C_0\bigoplus C_1 \to C_0 \to Y \to 0.$$ Again by Theorem
3.6, we get the following $\Hom_{\mathscr{A}}(\mathscr{C},-)$-exact
exact sequences:
$$C_0\bigoplus C_1\bigoplus C_2 \to C_0\bigoplus C_1 \to C_0 \to X \to 0$$
and
$$C_0\bigoplus C_1\bigoplus C_2 \to C_0 \to Y \to 0.$$
Continuing this procedure, we finally get the following
$\Hom_{\mathscr{A}}(\mathscr{C},-)$-exact exact sequences:
$$\cdots\to \bigoplus_{i=0}^{n-1}C_i\to \cdots \to C_0\bigoplus C_1\bigoplus C_2
\to C_0\bigoplus C_1 \to C_0 \to X \to 0$$
and
$$\cdots\to \bigoplus_{i=0}^{n-1}C_i\to \cdots \to C_0\bigoplus C_1\bigoplus C_2
\to C_0\bigoplus C_1 \to C_0 \to Y \to 0,$$ which implies that both
$X$ and $Y$ are objects in $\res\widetilde{\mathscr{C}}$.

(2) Let $M$ be an object in $\mathcal{G}(\mathscr{C})$ and
$$\cdots \to C_1\to C_0 \to C^0\to C^1\to \cdots$$ a
complete $\mathscr{C}$-resolution of $M$ with $M\cong\Im(C_0 \to
C^0)$. By (1) and Theorem 3.6, we get the following exact sequence:
$$\cdots\to \bigoplus_{i=0}^{n-1}C_i\to \cdots \to C_0\bigoplus C_1\bigoplus C_2
\to C_0\bigoplus C_1 \to C_0 \to X \to 0 \eqno{(4.1)}$$ which is
both $\Hom_{\mathscr{A}}(\mathscr{C},-)$-exact and
$\Hom_{\mathscr{A}}(-,\mathscr{C})$-exact. Dually, we get the
following exact sequence:
$$0\to X \to C^0 \to  C^0\bigoplus C^1\to C^0\bigoplus C^1\bigoplus C^2 \to
\bigoplus_{i=0}^{n-1}C^i \to \cdots \eqno{(4.2)}$$ which is also
both $\Hom_{\mathscr{A}}(\mathscr{C},-)$-exact and
$\Hom_{\mathscr{A}}(-,\mathscr{C})$-exact. Combining the sequences
(4.1) and (4.2), we conclude that $X$ is an object in
$\mathcal{G}(\mathscr{C})$. \hfill$\square$

\vspace{0.2cm}

We also have the following two out of three property.

\vspace{0.2cm}

{\bf Proposition 4.7.} {\it (1) If $\mathscr{C}$ is closed under
kernels of epimorphisms, then so is $\res\widetilde{\mathscr{C}}$.

(2) If $\mathscr{C}$ is closed under cokernels of monomorphisms,
then so is $\cores\widetilde{\mathscr{C}}$.

Let $$0\to X \to Y \to Z \to 0\eqno{(4.3)}$$ be an exact sequence in
$\mathscr{A}$.

(3) If the exact sequence (4.3) is
$\Hom_{\mathscr{A}}(\mathscr{C},-)$-exact and $X,Y$ are objects in
$\res\widetilde{\mathscr{C}}$, then $Z$ is also an object in
$\res\widetilde{\mathscr{C}}$.

(4) If the exact sequence (4.3) is
$\Hom_{\mathscr{A}}(-,\mathscr{C})$-exact and $Y,Z$ are objects in
$\cores\widetilde{\mathscr{C}}$, then $X$ is also an object in
$\cores\widetilde{\mathscr{C}}$.

(5) If the exact sequence (4.3) is both
$\Hom_{\mathscr{A}}(\mathscr{C},-)$-exact and
$\Hom_{\mathscr{A}}(-,\mathscr{C})$-exact, and if any two of $X,Y$ and
$Z$ are objects in $\mathcal{G}(\mathscr{C})$, then the third
term is also an object in $\mathcal{G}(\mathscr{C})$.}

\vspace{0.2cm}

{\it Proof.} We get (1) and (2) by Theorems 3.2 and
3.4 respectively. Dually, we get (3) and (4) by
Theorems 3.6 and 3.8 respectively.

(5) If  $X,Z$ are objects in $\mathcal{G}(\mathscr{C})$, then so is
$Y$ by [SSW1, Proposition 4.4]. By the above arguments, we get the
other two assertions. \hfill$\square$

\vspace{0.5cm}

\centerline{\large \bf 5. Gorenstein syzygies and dimension}

\vspace{0.2cm}

In this section, we fix subcategores $\mathscr{C}$ and $\mathscr{X}$
of $\mathscr{A}$. We use $\gen\mathscr{C}$ (resp.
$\cogen\mathscr{C}$) to denote a generator (resp. cogenerator)
for $\mathscr{C}$. The following result plays a crucial role in this
section.

\vspace{0.2cm}

{\bf Proposition 5.1.} {\it Let $\mathscr{C}$ be closed under
extensions and
$$0\to A\to C_1\buildrel {f} \over\to C_0\to M\to 0 \eqno{(5.1)}$$
an exact sequence in $\mathscr{A}$ with $C_0$ and $C_1$ objects in
$\mathscr{C}$.

(1) Then we have the following exact sequences:
$$0\to A\to C_1^{'} \to P_0\to M\to 0 \eqno{(5.2)}$$ and
$$0\to A\to I_1 \to C_0^{'}\to M\to 0 \eqno{(5.3)}$$
in $\mathscr{A}$ with $C_1^{'}, C_0^{'}$ objects in $\mathscr{C}$,
$P_0$ an object in $\gen\mathscr{C}$ and $I_1$ an object in
$\cogen\mathscr{C}$.

(2) For an object $D$ in $\mathscr{A}$, if (5.1) is
$\Hom_{\mathscr{A}}(D,-)$-exact (resp.
$\Hom_{\mathscr{A}}(-,D)$-exact), then so is (5.2) (resp. (5.3)).}

\vspace{0.2cm}

{\it Proof}. (1) There exists an exact sequence:
$$0\to C_0^{'} \to P_0 \to C_0\to 0$$ in $\mathscr{A}$ with $P_0$ an object in $\gen\mathscr{C}$
and $C_0^{'}$ an object in $\mathscr{C}$. Then we have the following
pull-back diagram:
$$\xymatrix{
& & 0 \ar[d] & 0 \ar[d] & & \\
& & C_0^{'} \ar[d] \ar@{=}[r] & C_0^{'} \ar[d] & \\
& 0 \ar[r] & N\ar[d] \ar[r] & P_0 \ar[d] \ar[r] & M \ar@{=}[d] \ar[r] & 0 \\
& 0 \ar[r] & \Im f\ar[d]\ar[r] & C_0 \ar[d] \ar[r] & M \ar[r] & 0 \\
& & 0 & 0 & & }$$ Consider the following pull-back diagram:
$$\xymatrix{
& & & 0 \ar[d]  & 0 \ar[d] &  \\
& & & C_0^{'} \ar[d] \ar@{=}[r] & C_0^{'} \ar[d] & \\
& 0 \ar[r] & A\ar@{=}[d] \ar[r] & C_1^{'} \ar[d] \ar[r] & N \ar[d] \ar[r] & 0 \\
& 0 \ar[r] & A\ar[r] & C_1 \ar[d]  \ar[r] & \Im f \ar[r]\ar[d] & 0 \\
& & & 0 &  0 & }$$ Because $\mathscr{C}$ is closed under extensions
and both $C_0^{'}$ and $C_1$ are objects in $\mathscr{C}$, $C_1^{'}$
is also an object in $\mathscr{C}$. Connecting the middle rows in
the above two diagrams, then we get the first desired exact
sequence. Dually, taking push-out, we get the second desired
sequence.

(2) If (5.1) is $\Hom_{\mathscr{A}}(D,-)$-exact, then so are the
middle rows in the above two diagrams by Lemma 2.4. Thus (5.2) is
also $\Hom_{\mathscr{A}}(D,-)$-exact. Dually, one gets the other
assertion. \hfill$\square$

\vspace{0.2cm}

{\bf Definition 5.2.}  Let $n$ be a positive integer. If there exists an exact sequence
$0\to A \to C_{n-1} \to C_{n-2} \to \cdots \to C_0 \to M \to 0$ in $\mathscr{A}$ with
all $C_i$ objects in $\mathscr{C}$, then $A$ is called an {\it $n$-$\mathscr{C}$-syzygy object}
(of $M$), and $M$ is called an {\it $n$-$\mathscr{C}$-cosyzygy object}
(of $A$).

\vspace{0.2cm}

The following theorem is one of main results in this section.

\vspace{0.2cm}

{\bf Theorem 5.3.} {\it Let $\mathscr{C}$ be closed under extensions,
and let $n\geq 1$ and
$$0\to A\to C_{n-1}\to C_{n-2}\to \cdots \to C_0\to M\to 0 \eqno{(5.4)}$$
be an exact sequence in $\mathscr{A}$ with all $C_i$ objects in
$\mathscr{C}$. Then we have the following

(1) There exist exact sequences:
$$0\to A\to I_{n-1}\to I_{n-2}\to\cdots\to I_0\to N\to 0 \eqno{(5.5)}$$ and
$$0\to M\to N\to X\to 0$$ in
$\mathscr{A}$ with all $I_i$ objects in $\cogen\mathscr{C}$ and $X$ an
object in $\mathscr{C}$. In particular, an object in $\mathscr{A}$
is an $n$-$\mathscr{C}$-syzygy if and only if it is an
$n$-$\cogen\mathscr{C}$-syzygy.

(2) There exist exact sequences:
$$0\to B\to P_{n-1}\to P_{n-2}\to\cdots\to P_0\to M\to 0 \eqno{(5.6)}$$
and $$0\to Y\to B\to A\to 0$$ in $\mathscr{A}$ with all $P_i$
objects in $\gen\mathscr{C}$ and $Y$ an object in $\mathscr{C}$. 
In particular, an object in $\mathscr{A}$
is an $n$-$\mathscr{C}$-cosyzygy if and only if it is an
$n$-$\gen\mathscr{C}$-cosyzygy.

(3) For an object $D$ in $\mathscr{A}$, if (5.4) is
$\Hom_{\mathscr{A}}(D,-)$-exact (resp.
$\Hom_{\mathscr{A}}(-,D)$-exact), then so are (5.5) (resp. (5.6)).}

\vspace{0.2cm}

{\it Proof.} (1) We proceed by induction on $n$.

Let $n=1$ and $$0\to A \to C_0\to M\to 0$$ be an exact sequence in
$\mathscr{A}$ with $C_0\in \mathscr{C}$. Then there exists an exact
sequence: $$0\to C_0 \to I_0\to X\to 0$$ in $\mathscr{A}$ with
$I_0$ an object in $\cogen\mathscr{C}$ and $X$ an object in $\mathscr{C}$.

Consider the following push-out diagram:
$$\xymatrix{
& & & 0 \ar[d] & 0 \ar[d] & \\
& 0 \ar[r] & A \ar@{=}[d] \ar[r] & C_0 \ar[d]\ar[r] & M \ar[d]\ar[r] & 0\\
& 0 \ar[r] & A\ar[r] & I_0 \ar[d]\ar[r] & N \ar[d]\ar[r] & 0 \\
& & & X \ar[d]\ar@{=}[r] & X \ar[d] & \\
& & & 0 & 0 & }$$ Then the middle row and the third column in the
above diagram are the desired exact sequences.

Now suppose that $n\geq 2$ and
$$0\to A\to C_{n-1}\to C_{n-2}\to \cdots \to C_0\to M\to 0$$
is an exact sequence in $\mathscr{A}$ with all $C_i$ objects in
$\mathscr{C}$. Put $K=\Coker(C_{n-1}\to C_{n-2})$. By Proposition 5.1(1), we
get an exact sequence:
$$0\to A\to I_{n-1}\to C_{n-2}^{'}\to K\to 0$$
in $\mathscr{A}$ with $I_{n-1}$ an object in $\cogen\mathscr{C}$ and
$C_{n-2}^{'}$ an object in $\mathscr{C}$. Put $A^{'}=\Im(I_{n-1}\to C_{n-2}^{'})$.
Then we get an exact sequence:
$$0\to A^{'}\to C_{n-2}^{'}\to C_{n-3} \to \cdots \to C_0\to M\to 0$$
in $\mathscr{A}$. Now we get the assertion by the induction hypothesis.

(2) The proof is dual to that of (1).

(3) It follows from Lemma 2.4 and Proposition 5.1 \hfill
$\square$

\vspace{0.2cm}

For any $n\geq 1$, we denote by $\Omega^n_{\mathscr{C}}(\mathscr{A})$
(resp. $\Omega^{-n}_{\mathscr{C}}(\mathscr{A})$) the subcategory of
$\mathscr{A}$ consisting of $n$-$\mathscr{C}$-syzygy (resp. $n$-$\mathscr{C}$-cosyzygy) objects.

\vspace{0.2cm}

{\bf Corollary 5.4.} {\it Let $\mathscr{C}$ be closed under
extensions. Then for any $n\geq 1$ we have

(1) If $\mathscr{X}$ is a cogenerator for
$\mathscr{C}$, then
$\Omega^n_{\mathscr{C}}(\mathscr{A})=\Omega^n_{\mathscr{X}}(\mathscr{A})$.

(2) If $\mathscr{X}$ is a generator for
$\mathscr{C}$, then
$\Omega^{-n}_{\mathscr{C}}(\mathscr{A})=\Omega^{-n}_{\mathscr{X}}(\mathscr{A})$.}

\vspace{0.2cm}

For an object $M$ in $\mathscr{A}$, the {\it
$\mathscr{C}$-dimension} of $M$, denoted by $\mathscr{C}$-$\dim M$,
is defined as inf$\{n\geq 0\mid$ there exists an exact sequence $0
\to C_{n} \to \cdots \to C_{1} \to C_{0} \to M \to 0$ in
$\mathscr{A}$ with all $C_i$ objects in $\mathscr{C}\}$. We set
$\mathscr{C}$-$\dim M$ infinity if no such integer exists.
A subcategory $\mathscr{X}$ of $\mathscr{C}$ is called
a {\it generator-cogenerator} for $\mathscr{C}$ if it is both a generator
and a cogenerator for $\mathscr{C}$.

Another main result in this section is the following

\vspace{0.2cm}

{\bf Theorem 5.5.} {\it Let $\mathscr{C}$ be closed under extensions
and $\mathscr{X}$ a generator-cogenerator for $\mathscr{C}$. Then
the following statements are equivalent for any object $M$ in
$\mathscr{A}$ and $n\geq 0$.

(1) $\mathscr{C}$-$\dim M \leq n$.

(2) There exists an exact sequence:
$$0\to C_n \to X_{n-1} \to \cdots \to X_1 \to X_0 \to M \to 0$$
in $\mathscr{A}$ with $C_n$ an object in $\mathscr{C}$ and all $X_i$
objects in $\mathscr{X}$.

(3) There exists an exact sequence:
$$0\to X_n \to X_{n-1} \to \cdots \to X_1 \to C_0 \to M \to 0$$
in $\mathscr{A}$ with $C_0$ an object in $\mathscr{C}$ and all $X_i$
objects in $\mathscr{X}$.

(4) For every non-negative integer $t$ such that $0 \leq t \leq n$,
there exists an exact sequence:
$$0\to X_{n}\to \cdots \to X_1\to X_0\to M\to 0$$
in $\mathscr{A}$ such that $X_t$ is an object in $\mathscr{C}$ and
all $X_i$ are objects in $\mathscr{X}$ for $i \neq t$.}

\vspace{0.2cm}

{\it Proof.} $(4)\Rightarrow (2)\Rightarrow (1)$ and $(4)\Rightarrow
(3)\Rightarrow (1)$ are trivial.

$(1)\Rightarrow (4)$ We proceed by induction on $n$.

Let $\mathscr{C}$-$\dim M \leq 1$ and $$0\to C_1\to C_0\to M\to 0$$
be an exact sequence in $\mathscr{A}$ with $C_0, C_1$ objects in
$\mathscr{C}$. By Proposition 5.1 with $A=0$, we get the exact
sequences $$0\to C_1^{'}\to X_0\to M\to 0$$ and $$0\to X_1\to
C_0^{'}\to M\to 0$$ in $\mathscr{A}$ with $X_0, X_1$ objects in
$\mathscr{X}$ and $C_0^{'}, C_1^{'}$ objects in $\mathscr{C}$.

Now suppose $\mathscr{C}$-$\dim M=n\geq 2$ and $$0\to C_{n}\to
\cdots \to C_1\to C_0\to M\to 0$$ is an exact sequence in
$\mathscr{A}$ with all $C_i$ objects in $\mathscr{C}$. Set
$A=\Coker(C_3\to C_2)$. By applying Proposition 5.1 to the exact
sequence: $$0\to A\to C_1\to C_0\to M\to 0,$$ we get the following
exact sequences:
$$0\to A \to C_1^{'}\to X_0\to M\to 0$$ and
$$0\to C_{n}\to \cdots \to C_2\to C_1^{'}\to X_0\to M\to 0$$
in $\mathscr{A}$ with $C_1^{'}$ an object in $\mathscr{C}$ and $X_0$
an object in $\mathscr{X}$. Set $N=\Coker (C_2\to C_1^{'})$. Then we
have $\mathscr{C}$-$\dim M \leq n-1$. By the induction hypothesis,
there exists an exact sequence: $$0\to X_{n}\to \cdots \to X_t\to
\cdots \to X_1\to X_0\to M\to 0$$ in $\mathscr{A}$ such that $X_t$
is an object in $\mathscr{C}$ and $X_i$ is an object in
$\mathscr{X}$ for $i\neq t$ and $1\leq t\leq n$.

In the following we only need to prove (4) for the case $t=0$. Set
$B=\Coker (C_2\to C_1)$. By the induction hypothesis, we get an
exact sequence: $$0\to X_n\to \cdots\to X_3 \to X_2\to C_1^{'} \to
B\to 0$$ in $\mathscr{A}$ with $C_1^{'}$ an object in $\mathscr{C}$
and all $X_i$ objects in $\mathscr{X}$. Set $K=\Coker (X_3\to X_2)$.
Then by applying Proposition 5.1 to the exact sequence: $$0\to K\to
C_1^{'}\to C_0\to M\to 0,$$ we get an exact sequence: $$0\to K\to
X_1\to C_0^{'}\to M\to 0$$ in $\mathscr{A}$ with $X_1$ an object in
$\mathscr{X}$ and $C_0^{'}$ an object in $\mathscr{C}$. Thus we
obtain the desired exact sequence: $$0\to X_n\to \cdots\to X_2 \to
X_1\to C_0^{'} \to M\to 0.$$ \hfill$\square$

\vspace{0.2cm}

Let $R$ be a ring. For a module $M$ in $\Mod R$ (resp. $\mod R$), we
use $\Add_RM$ (resp. $\add_RM$) to denote the subcategory of $\Mod R$
(resp. $\mod R$) consisting of all modules isomorphic to direct summands of
(finite) direct sums of copies of $_RM$. Recall that a subcategory
$\mathscr{D}$ of $\Mod R$ (resp. $\mod R$) is called {\it projectively
resolving} if $\mathscr{D}$ contains
(finitely generated) projective left $R$-modules, and $\mathscr{D}$
is closed under extensions and kernels of epimorphisms in
$\Mod R$ (resp. $\mod R$).

\vspace{0.2cm}

{\it Remark 5.6.} Let $R$ and $S$ be rings. Then we have

(1) Let $_RC_S$ be a semidualizing bimodule, that is, the following
conditions are satisfied: $(c1)$ $_{R}C$ admits a degreewise finite
$R$-projective resolution, and $C_{S}$ admits a degreewise finite
$S^{op}$-projective resolution; $(c2)$ both the homothety maps
$_{R}R_{R} \to \Hom_{S^{op}}(C, C)$
and $_{S}S_{S} \to \Hom_{R}(C, C)$
are isomorphisms; and $(c3)$ $\Ext_{R}^{i}(C,
C)=0=\Ext_{S^{op}}^{i}(C, C)$ for any $i \geq 1$. Then
$\mathcal{P}(\Mod R)\bigcup\Add_RC$ is a generator-cogenerator for
the subcategory of $\Mod R$ consisting of $G_C$-projective modules
(see [W] and [LHX, Corollary 2.10]).

(2) $\mathcal{P}(\Mod R)$ is a projective generator and an injective
cogenerator for the subcategory $\mathcal{G}(\Mod R)$ of Gorenstein
projective left $R$-modules. In particular, $\mathcal{G}(\Mod R)$ is
projectively resolving by [H, Theorem 2.5].

(3) If $R$ is an Artinian algebra and $T\in \mod R$ is cotilting,
then $\add_R(T\bigoplus R)$ is a generator-cogenerator
for $^{\bot}T$ by [AR1, Theorem 5.4(b)], and clearly $^{\bot}T$ is
projectively resolving, where $^{\bot}T=\{M\in \mod R\mid
\Ext_R^i(M,T)=0$ for any $i\geq 1\}$.

(4) Put $\Omega ^{n}(\Mod R)=\{M\in \Mod R\mid$ there exists an
exact sequence $0\to M \to P_{n-1} \to P_{n-2} \to \cdots \to P_0$ in $\Mod R$ with
all $P_i$ projective$\}$ for any $n\geq 1$,
and put $\Omega ^{\infty}(\Mod R)=\bigcap_{n\geq 1}\Omega ^{n}(\Mod R)$.
Then $\mathcal{P}(\Mod R)$ is a generator-cogenerator for
$\Omega ^{\infty}(\Mod R)$. Put $\Omega ^{n}(\mod R)$\linebreak $=\Omega ^{n}(\Mod R)\bigcap\mod R$
for any $n\geq 1$, and put $\Omega ^{\infty}(\mod R)=\bigcap_{n\geq 1}\Omega ^{n}(\mod R)$.
Then $\mathcal{P}(\mod R)$ is a
generator-cogenerator for $\Omega ^{\infty}(\mod R)$ over a left Noetherian ring $R$. By
[AR2, Theorem 1.7 and Proposition 2.2], we have that for a left and
right Noetherian ring $R$, if the right flat dimension of the $i$-th
term in a minimal injective coresolution of $R_R$ is at most $i+1$
for any $i\geq 0$, especially if $R$ is a commutative Gorenstein ring
(cf. [B, Fundamental Theorem]), then all $\Omega ^{n}(\mod R)$ and $\Omega ^{\infty}(\mod R)$
are closed under extensions.

(5) A subcategory $\mathscr{X}$ of $\mathscr{A}$ is called {\it
self-orthogonal}, denoted by $\mathscr{X}\bot\mathscr{X}$, if
$\Ext_{\mathscr{A}}^i(X,X^{'})=0$ for any objects $X$ and $X^{'}$ in
$\mathscr{X}$ and $i\geq 1$. If $\mathscr{X}\bot\mathscr{X}$, then
$\mathcal{G}(\mathscr{X})$ is closed under extensions and
$\mathscr{X}$ is a projective generator and an injective cogenerator
for $\mathcal{G}(\mathscr{X})$ by [SSW1, Corollaries 4.5 and 4.7].
In the rest of this section, we will focus on the self-orthogonal
subcategories of $\mathscr{A}$.

\vspace{0.2cm}

{\bf Lemma 5.7.} {\it Let $\mathscr{X}\bot\mathscr{X}$. Then
$\mathcal{G}(\mathscr{X})=(^{\bot}\mathscr{X}\bigcap\mathscr{X}^{\bot})\bigcap
(\res\widetilde{\mathscr{X}}\bigcap\cores\widetilde{\mathscr{X}})$.}

\vspace{0.2cm}

{\it Proof.} It is easy to get the assertion by the definition of
$\mathcal{G}(\mathscr{X})$. \hfill$\square$

\vspace{0.2cm}

As an application of Theorem 5.5, we have the following

\vspace{0.2cm}

{\bf Theorem 5.8.} {\it Let $\mathscr{X}\bot\mathscr{X}$ and $M$ be an
object in $\mathscr{A}$ with $\mathcal{G}(\mathscr{X})$-$\dim
M<\infty$. Then the following statements are equivalent for any
$n\geq 0$.

(1) $\mathcal{G}(\mathscr{X})$-$\dim M \leq n$.

(2) There exists an exact sequence:
$$0\to C_n \to X_{n-1} \to \cdots \to X_1 \to X_0 \to M \to 0$$
in $\mathscr{A}$ with $C_n$ an object in $\mathcal{G}(\mathscr{X})$
and all $X_i$ objects in $\mathscr{X}$.

(3) There exists an exact sequence:
$$0\to X_n \to X_{n-1} \to \cdots \to X_1 \to C_0 \to M \to 0$$
in $\mathscr{A}$ with $C_0$ an object in $\mathcal{G}(\mathscr{X})$
and all $X_i$ objects in $\mathscr{X}$.

(4) For every non-negative integer $t$ such that $0 \leq t \leq n$,
there exists an exact sequence:
$$0\to X_{n}\to \cdots \to X_1\to X_0\to M\to 0$$
in $\mathscr{A}$ such that $X_t$ is an object in
$\mathcal{G}(\mathscr{X})$ and all $X_i$ are objects in
$\mathscr{X}$ for $i \neq t$.

(5) $\Ext_{\mathscr{A}}^{n+i}(M,X)=0$ for any object $X$ in
$\mathscr{X}$ and $i\geq 1$.

(6) $\Ext_{\mathscr{A}}^{n+1}(M,Y)=0$ for any object $Y$ in
$\mathscr{A}$ with $\mathscr{X}$-$\dim Y<\infty$.

(7) $\Ext_{\mathscr{A}}^{n+i}(M,Y)=0$ for any object $Y$ in
$\mathscr{A}$ with $\mathscr{X}$-$\dim Y<\infty$ and $i\geq 1$.}

\vspace{0.2cm}

{\it Proof.} Because $\mathscr{X}\bot\mathscr{X}$ by assumption,
$\mathcal{G}(\mathscr{X})$ is closed under extensions and
$\mathscr{X}$ is a projective generator and an injective cogenerator for
$\mathcal{G}(\mathscr{X})$ by [SSW1, Corollaries 4.5 and 4.7]. So by
Theorem 5.5, we have $(1)\Leftrightarrow (2)\Leftrightarrow
(3)\Leftrightarrow (4)$.

$(7)\Rightarrow (6)$ and $(7)\Rightarrow (5)$ are trivial.

It is easy to get $(1)\Rightarrow (5)\Rightarrow (7)$ by Lemma 5.7
and the dimension shifting, respectively.

$(6)\Rightarrow (1)$ Let $\mathcal{G}(\mathscr{X})$-$\dim
M=m(<\infty)$. Then by Theorem 5.5, there exists an exact sequence:
$$0\to X_m \to X_{m-1} \to \cdots \to X_1 \to C_0 \to M \to 0$$
in $\mathscr{A}$ with all $X_i$ objects in $\mathscr{X}$ and $C_0$
an object in $\mathcal{G}(\mathscr{X})$. We claim that $m\leq n$.
Otherwise, let $m>n$. Note that $\mathscr{X}$-$\dim \Im(X_{n+1}\to
X_{n})\leq m-n-1<\infty$. So $\Ext_{\mathscr{A}}^1(\Im(X_n\to
X_{n-1}), \Im(X_{n+1}\to X_{n}))\cong \Ext_{\mathscr{A}}^{n+1}(M,
\Im(X_{n+1}\to X_{n}))=0$ (note: $X_0=C_0$) by assumption and Lemma
5.7. Hence the exact sequence $0\to \Im(X_{n+1}\to X_{n}) \to X_n
\to \Im(X_n\to X_{n-1}) \to 0$ splits and $\Im(X_n\to X_{n-1})$ is
isomorphic to a direct summand of $X_n$. By Theorem 4.6(2),
$\Im(X_n\to X_{n-1})$ is an object in $\mathcal{G}(\mathscr{X})$ and
$\mathcal{G}(\mathscr{X})$-$\dim M \leq n$, which is a
contradiction. \hfill$\square$

\vspace{0.2cm}

As an immediate consequence of Theorem 5.8, we get the following

\vspace{0.2cm}

{\bf Corollary 5.9.} {\it Let $\mathscr{X}\bot\mathscr{X}$ and $M$
an object in $\mathscr{A}$ with $\mathcal{G}(\mathscr{X})$-$\dim
M<\infty$. Then $\mathcal{G}(\mathscr{X})$-$\dim M=\sup\{n\geq 0\mid
\Ext_{\mathscr{A}}^n(M,X)\neq 0$ for some object $X$ in
$\mathscr{X}\}=\sup\{n\geq 0\mid \Ext_{\mathscr{A}}^n(M,Y)\neq 0$
for some object $Y$ in $\mathscr{A}$ with $\mathscr{X}$-$\dim
Y<\infty\}$.}

\vspace{0.2cm}

By Lemma 5.7 and Corollary 5.9, we get the following

\vspace{0.2cm}

{\bf Corollary 5.10.} {\it Let $\mathscr{X}\bot\mathscr{X}$ and
$0\to M_3\to M_2\to M_1 \to 0$ be an exact sequence in $\mathscr{A}$
with $M_3\neq 0$ and $M_1$ an object in $\mathcal{G}(\mathscr{X})$.
Then $\mathcal{G}(\mathscr{X})$-$\dim
M_3=\mathcal{G}(\mathscr{X})$-$\dim M_2$.}

\vspace{0.2cm}

In general, we have $\mathcal{G}(\mathscr{X})$-$\dim
M\leq\mathscr{X}$-$\dim M$ for any object $M$ in $\mathscr{A}$.
By Corollary 5.9 we get the following

\vspace{0.2cm}

{\bf Corollary 5.11.} {\it Let $\mathscr{X}\bot\mathscr{X}$ and $M$
an object in $\mathscr{A}$ with $\mathscr{X}$-$\dim M<\infty$. Then
$\mathcal{G}(\mathscr{X})$-$\dim M=\mathscr{X}$-$\dim M$.}

\vspace{0.2cm}

{\it Proof.} Let $\mathscr{X}$-$\dim M=n<\infty$. It suffices to
prove $\mathcal{G}(\mathscr{X})$-$\dim M\geq\mathscr{X}$-$\dim M=n$.
It is easy to get that $\Ext_{\mathscr{A}}^n(M,X)\neq 0$
some object $X$ in $\mathscr{X}$, so $\mathcal{G}(\mathscr{X})$-$\dim
M\geq n$ by Corollary 5.9. \hfill$\square$

\vspace{0.2cm}

By the definition of $\mathcal{G}(\mathscr{X})$, we have that each
object in $\mathcal{G}(\mathscr{X})$ can be embedded into an object
in $\mathscr{X}$ with the cokernel still in
$\mathcal{G}(\mathscr{X})$. The first assertion in the following
result generalizes this fact and [CFrH, Lemma 2.17].

\vspace{0.2cm}

{\bf Corollary 5.12.} {\it Let $\mathscr{X}\bot\mathscr{X}$ and $M$
be an object in $\mathscr{A}$ and $n\geq 0$, and let
$\mathcal{G}(\mathscr{X})$-$\dim M=n<\infty$. Then we have

(1) There exists an exact sequence $0\to M \to N \to G \to 0$ in
$\mathscr{A}$ with $\mathscr{X}$-$\dim N=n$ and $G$ an object in
$\mathcal{G}(\mathscr{X})$.

(2) There exists an exact sequence $0\to N\to G\to M\to 0$ in
$\mathscr{A}$ with $\mathscr{X}$-$\dim N\leq n-1$ and $G$ an object
in $\mathcal{G}(\mathscr{X})$. This exact sequence is an epic
$\mathcal{G}(\mathscr{X})$-precover of $M$.}

\vspace{0.2cm}

{\it Proof.} Let $M$ be an object in $\mathscr{A}$ with
$\mathcal{G}(\mathscr{X})$-$\dim M=n<\infty$.

(1) We apply Theorem 5.3(1) with $A=0$ to get an exact sequence
$0\to M\to N\to G\to 0$ in $\mathscr{A}$ with $\mathscr{X}$-$\dim
N\leq n$ and $G$ an object in $\mathcal{G}(\mathscr{X})$. By
Corollary 5.10, we have $\mathcal{G}(\mathscr{X})$-$\dim N=n$. Then
it follows from Corollary 5.11 that $\mathscr{X}$-$\dim N=n$.

(2) By Theorem 5.8, there exists an exact sequence $0\to N\to G\to
M\to 0$ in $\mathscr{A}$ with $\mathscr{X}$-$\dim N\leq n-1$ and $G$
an object in $\mathcal{G}(\mathscr{X})$. Also by Theorem 5.8,
$\Ext_{\mathscr{A}}^{1}(G^{'},N)=0$ for any object $G^{'}$ in
$\mathcal{G}(\mathscr{X})$. So the above exact sequence is an epic
$\mathcal{G}(\mathscr{X})$-precover of $M$. \hfill$\square$






\vspace{0.2cm}

{\it Remark 5.13.} For an object $M$ in $\mathscr{A}$, we may define
dually the {\it $\mathscr{C}$-codimension} of $M$, denoted by
$\mathscr{C}$-$\codim M$, as inf$\{n\geq 0\mid$ there exists an
exact sequence $0 \to M \to C^0 \to C^1 \to \cdots \to C^{n} \to 0$
in $\mathscr{A}$ with all $C^i$ objects in $\mathscr{C}\}$. We set
$\mathscr{C}$-$\codim M$ infinity if no such integer exists. We
point out the dual versions on the $\mathscr{C}$-codimension of all
the above results (5.5 and 5.8--5.12) also hold true by using a
completely dual arguments.

\vspace{0.5cm}

{\bf Acknowledgements.} This research was partially supported by the
Specialized Research Fund for the Doctoral Program of Higher
Education (Grant No. 20100091110034), NSFC (Grant Nos. 11171142,
11101217), NSF of Jiangsu Province of China (Grant Nos. BK2010047,
BK2010007) and a Project Funded by the Priority Academic Program
Development of Jiangsu Higher Education Institutions.

\end{document}